\documentclass[aos,seceqn,citesort,dvips]{arximspdf}

\doi{10.1214/10-AOS822}
\volume{39}
\issue{1}
\pubyear{2011}
\firstpage{439}
\lastpage{473}

\makeatletter

\newtheorem{theorem}{Theorem}[section]
\newtheorem{prop}[theorem]{Proposition}
\newtheorem{lemma}[theorem]{Lemma}
\newtheorem{cor}[theorem]{Corollary}
\newproclaim{example}{Example}[section]

\newcommand{\implies}{\Longrightarrow}

\newcommand{\mean}{\mathsf{E}}
\newcommand{\var}{\operatorname{\mathsf{Var}}}
\newcommand{\cov}{\operatorname{\mathsf{Cov}}}

\newcommand{\sr}{\mathsf{r}}
\newcommand{\sL}{\mathsf{L}}
\newcommand{\cH}{\mathcal{H}}
\newcommand{\cX}{\mathcal{X}}
\newcommand{\cZ}{\mathcal{Z}}

\newcommand{\tf}{\varphi}
\newcommand{\rx}{\varepsilon}
\newcommand{\px}{\vartheta}
\newcommand{\prx}{\theta}
\newcommand{\stx}{\sigma}
\newcommand{\state}{\eta}
\newcommand{\IK}{\kappa}

\newcommand{\Ints}{\mathbb{Z}}
\newcommand{\Nats}{\mathbb{N}}
\newcommand{\Reals}{\mathbb{R}}

\newcommand{\Gv}{ \mid}

\newcommand{\toi}{\to\infty}

\newcommand{\half}{\frac{1}{2}}

\newcommand{\cA}{\mathcal{A}}
\newcommand{\cB}{\mathcal{B}}

\makeatother

\begin{document}
\begin{frontmatter}

\title{Effects of statistical dependence on multiple testing under a
hidden Markov model}
\runtitle{Likelihood ratio under HMM}

\begin{aug}
\author[A]{\fnms{Zhiyi} \snm{Chi}\corref{}\thanksref{fund}\ead[label=e1]{zchi@stat.uocnn.edu}}
\runauthor{Z. Chi}
\affiliation{University of Connecticut}
\address[A]{Department of Statistics\\
University of Connecticut\\
215 Glenbrook Road, U-4120\\
Storrs, Connecticut 06269\\
USA\\
\printead{e1}} 
\end{aug}

\thankstext{fund}{Supported in part by NSF Grant DMS-07-06048 and NIH
MH 68028.}

\received{\smonth{4} \syear{2009}}
\revised{\smonth{3} \syear{2010}}

%
\begin{abstract}
The performance of multiple hypothesis testing is known to be
affected by the statistical dependence among random variables
involved. The mechanisms responsible for this, however, are not
well understood. We study the effects of the dependence structure
of a finite state hidden Markov model (HMM) on the likelihood
ratios critical for optimal multiple testing on the hidden states.
Various convergence results are obtained for the likelihood ratios
as the observations of the HMM form an increasing long chain.
Analytic expansions of the first and second order derivatives are
obtained for the case of binary states, explicitly showing the
effects of the parameters of the HMM on the likelihood ratios.
\end{abstract}

%

\end{frontmatter}

\section{Introduction} \label{sec:intro}
Statistical dependence in data poses a challenge to multiple
hypothesis testing. Under the framework of the false discovery rate
(FDR), many efforts have been made to establish the control of FDR
under dependence \cite
{benjaminiyek01,farcomeni07,wufdr08,suncai09,sarkar06}.
Meanwhile, many empirical and analytical
works have described the effects of dependence on the outputs of
multiple tests \cite{owen05,qiu05,efron07b,halljin08}.
However, in what way the dependence impacts multiple testing is not
well understood.

A useful model that incorporates tractable dependence in multiple
testing is the hidden Markov model (HMM) \cite{suncai09}. In the
model, the nulls are organized as $H_t$, where the index $t$ takes
integer values. Each $H_t$ is associated with a random variable
that determines whether the null is true or false. The random
variables form a Markov chain but are hidden and unobservable.
Instead, the observations $X_t$ each is a many-to-one transform of
the hidden variable corresponding to $H_t$. In the context of
multiple testing, it will be useful to treat the hidden variable as
consisting of two parts, $\state_t$ and $Z_t$. On the one hand,
$\state_t$ encodes the ``true identity,'' or state of the signal
associated with $H_t$ and in general can take two or more possible
values. On the other, $Z_t$ acts as the noise that blurs or distorts
the signal. Then $X_t$ can be thought of as the result of a
deterministic interaction between $\state_t$ and $Z_t$.

To understand the role of dependence in the multiple tests on the
nulls, the ``oracle'' approach assumes the parameters in the HMM are
known and explores what amounts to an optimal testing procedure. The
advantage of this approach is that it can reveal effects purely due to
dependence, without confounding with effects due to specific parameter
estimation methods. Suppose the observations are $X_{-m}, \ldots,
X_n$. With the parameters being known, for each null $H_t$, the
conditional likelihood $\operatorname{\mathsf{Pr}}\{\mbox{$H_t$ is true}\mid X_{-m},
\ldots, X_n\}$ can be computed. The importance of the conditional likelihood
for multiple testing has been shown in various contexts
\cite
{efronetal01,mulleretal04,storey07opt,chi08mvp,suncai09}. For the
HMM, \cite{suncai09} shows that under a certain loss function, an
optimal procedure is to reject $H_t$ if and only if the corresponding
conditional likelihood is small enough. The loss function is a linear
combination of the numbers of Types I and II errors and is related to
the FDR. The importance of the conditional likelihood can also be
argued directly based on the FDR criterion, and in fact without
particular assumption on dependence; see the \hyperref[app]{Appendix}.

In view of the role of the conditional likelihood, our aim is to
investigate how it is affected by the parameters of the HMM. The
parameters can be divided into two types, respectively, characterizing
the dependence among $\state_t$ and the ``strength'' of useful
signals. In addition, the conditional likelihood also depends on how
$\state_t$ and $Z_t$ interact. The next example illustrates what role
may be expected for these factors.
\begin{example}\label{ex:gauss}
Suppose the states $\state_t$ are equal to $\mathbf{1}\{H_t \mbox{
is false}\}$ and form a stationary Markov chain with transition
probabilities $q_{ij} = \mathsf{Pr}\{\state_t=j\mid\state_{t-1}=i\}
>0$; moreover,
conditional on $\state=(\state_t)$, $X_t$ are independent $\sim
N(\rx\state_t,1)$ with $\rx>0$. Write $X_t=Z_t + \rx\state_t$.
Then $(Z_t, \state_t)$ form a hidden Markov chain, with $Z_t$ i.i.d.
$\sim N(0,1)$. The strength of the signals is measured by $\rx$,
the interaction between the noise $Z_t$ and $\state_t$ is additive,
such that $X_t = \tf(Z_t, \rx\state_t)$ with $\tf(z,\px) = z+\px$.

In many cases, the observations form a long chain $X_{-m}, \ldots,
X_n$, with $m$, $n\gg1$, so the effect of the parameters can be
studied through the properties of
\[
\mathsf{Pr}\{\state_t=0\mid X\}
= \lim_{m,n\toi} \operatorname{\mathsf{Pr}}\{\state_t=0\mid X_{-m}, \ldots, X_n\}
\]
for each $t$, where $X=(X_t, t\in\Ints)$. Since
$\mathsf{Pr}\{\state_t=0\mid X_{-m}, \ldots, X_n\}$ form a martingale
for any
increasing sequence of $m$ and $n$, the (almost sure) existence of
the limit is guaranteed. However, this says nothing about how the
limit depends on $\rx$ and $q_{ij}$. To get some insight, consider
instead the likelihood ratios
\[
\frac{\mathsf{Pr}\{\state_t=1\mid X\}}{\mathsf{Pr}\{\state_t=0\mid X\}}
= \frac{1}{\mathsf{Pr}\{\state_t=0\mid X\}}-1,
\]
which turn out to be a little more convenient to study.
Regarding them as functions of $\rx$, we next
consider their Taylor expansions. In principle, the likelihood
ratios can be expanded around any value of $\rx$.
Since large values of $|\rx|$ correspond to strong signals whose
detection is easy, we shall expand around $\rx=0$ to get
insight into the case where the strength of signal ranges from
being weak\vadjust{\goodbreak} to moderate. Without loss of
generality, consider the likelihood ratio for $\state_0$. Since
$\state$ is stationary, $\mathsf{Pr}\{\state_{t-1}=j\mid\state_t=i\}=
q_{ij}$. By
the Bayes rule and Markov property,
\begin{eqnarray*}
&&\mathsf{Pr}\{\state_0=a\mid X_{-m}, \ldots, X_n\}\\
&&\qquad\propto
P(a)\mathop{\sum_{\stx_{-m}, \ldots, \stx_n}}_{\stx_0=a}
\exp\Biggl\{-\frac{1}{2}\sum_{t=-m}^n (Z_t + \rx\state_t - \rx\stx
_t)^2\Biggr\}
\prod_{t=0}^{n-1} q_{\stx_t\stx_{t+1}}
\prod_{t=0}^{m-1} q_{\stx_{-t}\stx_{-t-1}}
\end{eqnarray*}
for $a=0,1$, where $P(a)=\mathsf{Pr}\{\state_0=a\}$. Then, formally, one
can get
\begin{eqnarray*}
\frac{d}{d\rx}
\biggl[\ln\frac{\mathsf{Pr}\{\state_0=1\mid X\}}{\mathsf{Pr}\{\state_0=0
\mid X\}} \biggr]_{\rx=0}
&=&
\lim_{m,n\toi}
\frac{d}{d\rx}
\biggl[\ln\frac{\mathsf{Pr}\{\state_0=1\mid X_{-m}, \ldots, X_n\}}
{\mathsf{Pr}\{\state_0=0\mid X_{-m}, \ldots, X_n\}} \biggr]_{\rx=0} \\
&=&
\sum_{t=-\infty}^\infty
Z_t[\mathsf{Pr}\{\state_t=1\mid\state_0=1\} -\mathsf{Pr}\{\state
_t=1\mid\state_0=0\} ] \\
&=& \sum_{t=-\infty}^\infty r^{|t|} Z_t,
\end{eqnarray*}
where $r\not=1$ is one of the two eigenvalues of the matrix $(q_{ij})$,
the other being 1.

We shall refer to the above conditional likelihood ratio as the full
likelihood ratio (FLR), as it is based on the entire $X$. On the
other hand, if the information of the dependence (i.e., $q_{ij}$) is
not available, but the values of all other parameters are known,
including $P(a)$, then the likelihood ratio would have to be
evaluated as
\[
\frac{\mathsf{Pr}\{\state_0=1\Gv X_0\}}{\mathsf{Pr}\{\state_0=0\Gv
X_0\}}
=
\frac{P(1)}{P(0)} \frac{f(X_0-\rx)}{f(X_0)}=
\frac{P(1)}{P(0)} \exp\biggl\{\rx Z_0 + \frac{\rx^2}{2}(2\state
_0-1)\biggr\},
\]
where $f(x)$ is the density of $N(0,1)$. We shall refer to this
conditional likelihood ratio as the local likelihood ratio (LLR).
It then can be seen that
\[
\ln\frac{\mathrm{FLR}}{\mathrm{LLR}}
= \rx\sum_{t\not=0} r^{|t|} Z_t + o(\rx).
\]
Thus, at the first order, the dependence in $\state$ merely adds
noise but no ``net effect,'' regardless of the actual values of
$\state$. If there is any state-dependent effect, it
should be reflected in a higher order term of $\rx$. To see if this
is the case, take the second order derivative in $\rx$. Again, the
calculation can be done formally. To evaluate the state-dependent
net effect, proceed with
\begin{eqnarray*}
\mean[(\ln\mathrm{FLR})''_{\rx=0} \Gv\state_0 ]
&=&
\lim_{m,n\toi}
\mean\biggl[\frac{d^2}{d\rx^2} \biggl[\ln\frac{\mathsf{Pr}\{\state
_0=1\mid X_{-m}, \ldots, X_n\}} {\mathsf{Pr}\{\state_0=0\mid X_{-m},
\ldots,
X_n\}} \biggr]_{\rx=0}\Bigm|\state_0 \biggr] \\
&=& (2\state_0-1) \sum_t r^{|2t|},
\end{eqnarray*}
giving
\[
\mean\biggl[\biggl(\ln\frac{\mathrm{FLR}}{\mathrm{LLR}} \biggr)''_{\rx=0}\Bigm|\state_0
\biggr] = (2\state_0-1) \sum_{t\not=0} r^{2|t|}.
\]
It follows that, comparing to $\ln\mathrm{LLR}$, if $\state_0=1$, on
average $\ln\mathrm{FLR}$ is larger, making $H_0$ more likely to be
(correctly) rejected, whereas if $\state_0=0$, it is smaller, making
$H_0$ less likely to be (falsely) rejected.

So far the expansions are expressed in terms of the unobservable
$Z_t$. One question is whether similar expansions in terms of
the observable $X_t$ can be obtained. As will be seen in Section
\ref{ssec:univariate}, this is possible after we get more
information on higher order derivatives.

From the expansions, the effect of the dependence in $\state$ on the
likelihood ratio is apparent. In both the first and second order
derivatives, the effect is determined by~$r$. In particular, when
$r=0$, $\state_t$ are i.i.d. and FLR is equal to LLR. Consistent with
this, the derivatives of the difference between the two ratios
become 0.
\end{example}

As the example, the rest of the paper develops Taylor's expansion in
terms of $\rx$ for the FLRs $\mathsf{Pr}\{H_t \mbox{ is false}\mid X\}
/\mathsf{Pr}\{H_t \mbox{ is true}\mid X\}$ to study the effects of
dependence structure of HMM.
The differentiation involved in the expansion
should be interpreted as follows. During the differentiation,
both the signal $\state$ and noise $Z$ are fixed. As the strength
$\rx$ of the signal varies, the observed values $X_t$ become functions
of $\rx$. The likelihood ratio is affected by $\rx$ in two ways:
not only the value of $X_t$ is changed, but also the parametric form
of the density function of $X_t$. Both have to be taken into
account in the derivatives.

Several issues need to be addressed. First, we have only considered a
stationary process of the signals $\state$. In applications, it is
useful to consider nonstationary $\state$ that has time-dependent
transition probabilities. Moreover, it is useful to consider
various types of interactions between $\state_t$ and $Z_t$ besides the
additive one.

Second, in Example \ref{ex:gauss}, each $\state_t$ is binary,
indicating whether a null is true or false. For more generality, one
can assume a finite state Markov chain, such that a subset
of the states are associated with true nulls and the rest with false
nulls. Even for a binary process, it can be useful
to
reformulate it as a multistate Markov chain. For example, let
$\state$ be a second order binary Markov chain, that is, $\mathsf
{Pr}\{\state_t\mid\state_s, s<t\} = \mathsf{Pr}\{\state_t\mid\state
_{t-1}, \state_{t-2}\}$. Then one can define a first order Markov chain
$\tilde\state$ by $\tilde\state_t = (\state_{t-1},
\state_t)$. If $\state_t = \mathbf{1}\{H_t \mbox{ is false}\}$,
then in
terms of $\tilde\state$, $(0,0)$ and $(1,0)$ are states associated
with true nulls, and $(0,1)$ and $(1,1)$ are states associated with
false nulls.

Third, in Example \ref{ex:gauss}, limit operation, differentiation,
and expectation are freely interchanged for $\mathsf{Pr}\{\state_t\mid
X_{-m}, \ldots, X_n\}$ for fixed $t$. This has to be justified. Note that
the likelihood bears similarity to $\mathsf{Pr}\{\state_n\mid X_0,
\ldots, X_n\}$, a quantity extensively studied in the literature on
nonlinear filtering
and related issues \cite
{atar97,atar97b,baxendale04,chigansky06,chigansky04,dimasi05,douc04,douc01,fleming99,kaijser75,kochman06,legland00,legland00b,tadic05}.
As in
most of the cited works, in this paper, convergence results are
established using geometric contraction. On the other hand, in
those works, the goal is to establish weak convergence of the
conditional probability for $\state_n$ under the assumption of
stationary transition probabilities. As seen in Example
\ref{ex:gauss}, the convergence of the conditional
probability for $\state_t$ follows from the martingale
convergence. So instead, the goal here is to establish convergence
for the derivatives of the conditional likelihood with
arbitrary transition probabilities.

The rest of the paper proceeds as follows. In Section \ref{sec:main},
a HMM is set up in the context of multiple testing and then
various convergence results on the likelihood ratio are stated. In
Section \ref{sec:example}, the likelihood ratio for a first order HMM
with binary states is considered in more detail, which allows more
explicit expressions for the first and second derivatives of the
likelihood ratio. Several examples are given, with Example
\ref{ex:gauss} being a special case. Theoretical details are provided
in Section~\ref{sec:proof}.

\section{Main results} \label{sec:main}
\subsection{A HMM setup}
Let $\state=\{\state_t, t\in\Ints\}$ be a finite state process,
such that the state space $\cH$ is partitioned into $\cH_0$ and
$\cH_1$, with states in $\cH_0$ being associated with true nulls,
while those in $\cH_1$ associated with false nulls. The noise process
is $Z=\{Z_t, t\in\Ints\}$, with each $Z_t$ taking values in a
Euclidean space $\cZ$. To model the interaction between $\state_t$
and $Z_t$, let $\{\tf(z, \px), \px\in\Theta\}$ be a family of
mappings $\cZ\to\cX$ indexed by $\px$, where $\Theta$ is an open set
in $\Reals^d$ and $\cX$ a Euclidean space. Then, let
\[
\prx_a\dvtx\Reals^p \to\Theta,\qquad a\in\cH,
\]
be a family of functions, such that each $\rx\in\Reals^p$ specifies a
scenario where the observations are
%
%
\begin{equation} \label{eq:observe}
X_t=X_t(\rx) = \tf(Z_t, \prx_{\state_t}(\rx)).
\end{equation}

Intuitively, $\tf(Z_t, \px)$ determines how $Z_t$ interacts with a
possible manifestation of $\state_t$ to generate an observation $X_t$;
the manifestation of $\state_t$ is $\prx_{\state_t}(\rx)$, with
$\rx$
being the tuning parameter that determines how strongly $\state_t$
manifests itself. The dimension $p$ of $\rx$ may be greater than 1 to
take into account different aspects of the tuning. We will assume
that $(\state,Z)$ is defined on the canonical space
$\cH^\Ints\times\cZ^\Ints$ equipped with the product Borel
$\sigma$-algebra.

For function $h\dvtx\Reals^s\to\Reals$ and $s$-tuple of nonnegative
integers $\nu=(\nu_1, \ldots, \nu_s)$, denote the $\nu$th
derivative of $h$ and
its order, respectively, by
\[
h^{(\nu)}(x) = \frac{\partial^{|\nu|} h(x)}{\partial x_1^{\nu_1}
\cdots\partial x_s^{\nu_s}},\qquad |\nu|=\nu_1+\cdots+\nu_s.
\]
Denote $h^{(\nu)}:=h$ if $\nu=0:=(0,\ldots,0)$. For $q\in\Nats$,
denote $h\in C^{(q)}$ if $h^{(\nu)}$ exists and is continuous for every
$|\nu|\le q$. If $i=(i_1, \ldots, i_{s})$ and $\nu=(\nu_1, \ldots, \nu_
s)$, denote $i\le
\nu$ if $i_k\le\nu_k$ for every $k=1,\ldots,s$ and denote $i<\nu$ if
$i\le\nu$ and $i\not=\nu$.

\subsubsection*{Assumptions}
The following assumptions will be needed for different occasions:
\begin{enumerate}
\item\hypertarget{A:noise}
$Z$ is independent of $\state$ and $Z_t$ are i.i.d. such that for each
$\px\in\Theta$ and $t\in\Ints$, $\tf(Z_t,\px)$ has a
density $f(x,\px)$.
\item\hypertarget{A:state}
$\state$ is a Markov chain and there are $\IK\ge1$, $\phi_*>0$,
such that for all $a$, $b\in\cH$ and $s$, $t\in\Ints$ with
$|s-t|\ge\IK$, $\mathsf{Pr}\{\state_t=b\mid\state_s=a\} \ge\phi_*$.
\item\hypertarget{A:positive}
For each $z\in\cZ$ and $a$, $b\in\cH$, $0<f(\tf(z, \prx_a(\rx)),
\prx_b(\rx))<\infty$ and is continuous in $\rx$.
\item\hypertarget{A:density}
There is $q\in\Nats$, such that for each $z\in\cZ$ and $a$,
$b\in\cH$, $f(\tf(z, \prx_a(\rx)), \prx_b(\rx))$ as a function in
$\rx$ belongs to $C^{(q)}$ and all its partial derivatives of order
$\le q$ are continuous in $(z,\rx)$. Furthermore, for $r>0$,
there is $c=c(r)>2$, such that
\[
\mathsf{Pr}\{M_q(Z_0, r) \ge u\} = O((\log u)^{-c}),\qquad u\toi,
\]
where, letting $\ell_{z,ab}(\rx) = \ln f(\tf(z, \prx_a(\rx)),
\prx_b(\rx))$, $M_0(z,r)=1$ and for $k>0$,
\[
M_k(z,r) = \sup\bigl\{\bigl|\ell_{z,a b}^{(\nu)}(\rx)\bigr|\dvtx1\le|\nu|\le k,
|\rx|\le r, a,b\in\cH\bigr\}.
\]
\item\hypertarget{A:density2}
For any $r>0$, $E[M_q(Z_0, r)^k]<\infty$, where $k=q^2(q+1)/2$.
\end{enumerate}
Henceforth, for $s$, $t\in\Ints$ and $a$, $b\in\cH$, denote
\[
P_t(a) = \mathsf{Pr}\{\state_t=a\},\qquad
P_{st}(a,b) = \mathsf{Pr}\{\state_t=b\mid\state_s=a\}.
\]

\subsubsection*{Remarks}
1. Some examples of $\tf$ are given Section \ref{ssec:examples}.

2. $\eta$ need not be stationary or have time-homogeneous
transitions.

3. Assumption \hyperlink{A:positive}{3} implies that no value of $X_t$ can
decisively identify or rule out any elements in $\cH$ as possible
values for $\state_t$.

4. In Example \ref{ex:gauss}, since $\ell_{z,a b}(\rx)= -\frac{1}{2}
[z+\rx(a-b)]^2 - \ln\sqrt{2\pi}$ and $Z_t\sim N(0,1)$,
Assumption \hyperlink{A:density2}{5} is satisfied. The assumption is stronger
than Assumption \hyperlink{A:density}{4}. To get results on almost sure
convergence, Assumption \hyperlink{A:density}{4} suffices. However, to get
results on expectations, Assumption \hyperlink{A:density2}{5} will be used.

5. Assumption \hyperlink{A:state}{2} can be relaxed as follows: there are
$\phi_*>0$ and $\cdots<s_k<t_k < s_{k+1}<t_{k+1}<\cdots,$ with $s_k \to
\pm\infty$ as $k\to\pm\infty$, such that $P_{s_k,t_k}(a,b)\ge\phi_*$
and for $n\gg1$, $\#\{k\dvtx-n\le s_k\le0\}/n$ and $\#\{k\dvtx0\le
s_k\le
n\}/n$ are bounded away from 0. The analysis under the relaxed
assumption follows the same line as the rest of the paper but is more
technical. We will not pursue it here.

\subsection{Derivatives of full likelihood ratios}
Given $\rx$ and $m$, $n\in\Nats$, if the observations consist
of $X_s(\rx) = \tf_s(Z_s, \prx_{\state_s}(\rx))$ with $s=-m,
\ldots,
n$, the likelihood ratio for false null vs. true null at
$t$ is
\[
\rho_{t,m n}(\rx)
=
\frac{
\mathsf{Pr}\{\state_t\in\cH_1\mid X_{-m}(\rx), \ldots, X_n(\rx)\}
}{
\mathsf{Pr}\{\state_t\in\cH_0\mid X_{-m}(\rx), \ldots, X_n(\rx)\}
}.
\]

Let $\stx=(\stx_t)$ be an independent copy of $\state$ that is
independent of $Z$ as well. Denote by $\mean_\sigma$ the expectation
with respect to $\stx$. By Bayes formula,
%
%
\begin{equation} \label{eq:LR}
\rho_{t,m n}(\rx)
= \frac{
\sum_{a\in\cH_1} P_t(a)
\mean_\sigma[\prod_{s=-m}^n \psi_s(\rx,\stx_s)\mid\stx_t=a]
}{
\sum_{a\in\cH_0} P_t(a)
\mean_\sigma[\prod_{s=-m}^n \psi_s(\rx,\stx_s)\mid\stx_t=a]
},
\end{equation}
where for $c\in\cH$,
%
%
\begin{equation} \label{eq:psi-t}
\psi_t(\rx,c) = f(X_t(\rx), \prx_c(\rx)) =
f(\tf(Z_t, \prx_{\state_t}(\rx)), \prx_c(\rx)).
\end{equation}

As discussed in the \hyperref[sec:intro]{Introduction},
\[
\rho_t(\rx)=\lim_{m,n\toi}\rho_{t,m n}(\rx) =
\frac{\mathsf{Pr}\{\state_t\in\cH_1\mid X_s(\rx), s\in\Ints\}}
{\mathsf{Pr}\{\state_t\in\cH_0\mid X_s(\rx), s\in\Ints\}}
\]
exists almost surely due to martingale convergence and plays an
important role in optimal multiple testing procedures.
\begin{theorem} \label{thm:LR-limit}
Suppose Assumptions \hyperlink{A:noise}{1}--\hyperlink{A:density}{4} hold.

1. Almost surely, $\rho_{t,m n}\in C^{(q)}$ for $t=-m+\IK,
\ldots, n-\IK$.

2. Almost surely, $\rho_t(\rx)$ is strictly positive for all $t$
and $\rx$.

3. There is a deterministic function $r_{t,\nu}(\rx_0)\in
(0,1)$ in $\rx_0>0$ for each $t\in\Ints$ and $\nu$ with $|\nu|\le
q$, such that almost surely, as $m$, $n\toi$, $\rho_{t,m n}^{(\nu
)}(\rx)$ converges, with
\[
\sup_{|\rx|\le\rx_0} \Bigl|\rho_{t,m n}^{(\nu)}(\rx) - \lim_{m,n\toi
} \rho_{t,m n}^{(\nu)} (\rx)\Bigr|= o(r_{t,\nu}^{m\wedge n}(\rx_0)),
\]
for all $t\in\Ints$, $\nu$ with $|\nu|\le q$ and $\rx_0>0$.
\end{theorem}

Due to the uniform convergence of $\rho_{t,m n}^{(\nu)}$ on
every compact set,
%
%
\begin{equation} \label{eq:uniform}
\rho_t\in C^{(q)},\qquad \rho_t^{(\nu)}(\rx) = \lim_{m,n\toi}
\rho_{t,m n}^{(\nu)}(\rx),\qquad t\in\Ints, |\nu|\le q
\end{equation}
(cf. \cite{rudin76}, Theorem 7.17). Then, as $\rho_t(\rx)$ are
strictly positive, the interchange between limit operation and
differentiation for the logarithms of
$\rho_{t,m n}(\rx)$
in Example~\ref{ex:gauss} is
justified.

Since $\Ints$ is countable, in order to establish Theorem
\ref{thm:LR-limit}, it suffices to show it holds for each fixed
$t\in\Ints$. Without loss of generality, we shall focus on $t=0$.
For ease of notation, henceforth denote $\rho_{m n}=\rho_{0,m n}$.

By the conditional independence of $(\stx_t, t<0)$ and $(\stx_t, t>0)$
given $\stx_0$,
\begin{eqnarray*}
\mean_\stx\Biggl[\prod_{s=-m}^n \psi_s(\rx,\stx_s)\biggm|\stx_0\Biggr]
&=&
\psi_0(\rx, \stx_0)
\mean_\stx\Biggl[\prod_{s=1}^n \psi_s(\rx,\stx_s)\biggm|\stx_0\Biggr]\\
&&{}\times
\mean_\stx\Biggl[\prod_{s=1}^m \psi_{-s}(\rx,\stx_{-s})\biggm|\stx_0\Biggr].
\end{eqnarray*}

Fix an arbitrary $\imath\in\cH$. Define
%
%
\begin{equation} \label{eq:Lambda}
\Lambda_{\pm n,a} = \Lambda_{\pm n,a}(\rx) =
\frac{
\mean_\stx[\prod_{s=1}^n \psi_s(\rx,\stx_{\pm s})\mid\stx_0=a]}
{\mean_\stx[\prod_{s=1}^n \psi_s(\rx,\stx_{\pm s})\mid\stx
_0=\imath]}.
\end{equation}
Then (\ref{eq:LR}) for $t=0$ can be written as
%
%
\begin{equation} \label{eq:rhom n}
\rho_{m n}(\rx)=
\frac{
\sum_{a\in\cH_1} \psi_0(\rx,a) P_0(a)
\Lambda_{-m,a}\Lambda_{n,a}
}{
\sum_{a\in\cH_0} \psi_0(\rx,a) P_0(a)
\Lambda_{-m,a}\Lambda_{n,a}
}.
\end{equation}
%
Although $\Lambda_{\pm n, a}$ depends on $\imath$, $\rho_{mn}(\rx)$ is independent of $\imath$.
For brevity, $\imath$ is omitted in the notation.

From (\ref{eq:rhom n}), it is seen that Theorem \ref{thm:LR-limit}
follows from the next two assertions on uniform geometric contraction
of functions and their derivatives on any compact interval of $\rx$.
\begin{theorem} \label{thm:L-limit}
Let Assumptions \hyperlink{A:noise}{1}--\hyperlink{A:positive}{3} hold. Almost
surely, as $n\toi$, for all $a\in\cH$, $\Lambda_{\pm n,a}(\rx)$
converge uniformly on every compact set
of $\rx$. The limits
%
%
\begin{equation} \label{eq:Lambda-lim}
\sL_a(\rx) = \lim_{n\toi} \Lambda_{n,a}(\rx),\qquad
\bar\sL_a(\rx) = \lim_{n\toi} \Lambda_{-n,a}(\rx)
\end{equation}
are strictly positive and continuous, and there is a deterministic
increasing function $r(\rx_0)\in(0,1)$ in $\rx_0>0$, such that
almost surely, as $n\toi$,
\[
{\sup_{|\rx|\le\rx_0}}|\Lambda_{n,a}(\rx) - \sL_a(\rx) |=o(r(\rx
_0)^{n})\qquad \forall\rx_0>0,
\]
and likewise for $\Lambda_{-n,a}$ and $\bar\sL_a(\rx)$.
\end{theorem}
\begin{theorem} \label{thm:d-L-limit}
Let Assumptions \hyperlink{A:noise}{1}--\hyperlink{A:density}{4} hold. Almost
surely, as $n\toi$, for each nonzero $\nu$ with $|\nu|\le q$ and
$a\in\cH$, $\Lambda_{\pm n,a}^{(\nu)}$ converge uniformly on every
compact set of $\rx$. Let
\[
\sL_{\nu,a}(\rx) = \lim_{n\toi} \Lambda_{n,a}^{(\nu)}(\rx),\qquad
\bar\sL_{\nu,a}(\rx) = \lim_{n\toi} \Lambda_{-n,a}^{(\nu)}(\rx).
\]
There is an increasing deterministic function $r_\nu(\rx_0)\in
(0,1)$ in $\rx_0>0$, such that almost surely, as $n\toi$,
\[
\max_a\sup_{|\rx|\le\rx_0}\bigl|\Lambda_{n,a}^{(\nu)}(\rx) - \sL
_{\nu,a}(\rx) \bigr|=o(r_\nu^n(\rx_0))\qquad \forall\rx_0>0,
\]
and likewise for $\Lambda_{-n,a}$ and $\bar\sL_{\nu,a}(\rx)$.
\end{theorem}

Basically, the two theorems say that $\sL_a(\rx)$ and
$\bar\sL_a(\rx)$ are $q$ times differentiable, and for $\nu$ with
$|\nu|\le q$, $\sL_a^{(\nu)}(\rx) = \sL_{\nu,a}(\rx)$, $\bar\sL
_a^{(\nu)}
(\rx) = \bar\sL_{\nu,a}(\rx)$, that is, $(\lim\Lambda_{\pm
n,a})^{(\nu)}= \lim\Lambda_{\pm n,a} ^{(\nu)}$. As a result,
$\rho(\rx)$ is $q$ times differentiable, with
%
%
\begin{equation} \label{eq:rho-limit}
\rho^{(\nu)}(\rx) =
\biggl[\frac{ \sum_{a\in\cH_1} \psi_0(\rx,a) P_0(a) \sL_a(\rx)\bar
\sL_a(\rx) }{ \sum_{a\in\cH_0} \psi_0(\rx,a) P_0(a) \sL_a(\rx
)\bar\sL_a(\rx) } \biggr]^{(\nu)}.
\end{equation}
Note that although we are mainly interested on the property of
$\rho_t$ around $\rx=0$, the above results allow Taylor's expansion
around nonzero values of $\rx$ as well.\vadjust{\goodbreak}

In Example \ref{ex:gauss}, limit operation, differentiation, and
expectation were freely interchanged. The next assertion justifies
this.
\begin{theorem}\label{thm:k=1}
Let Assumptions \hyperlink{A:noise}{1}--\hyperlink{A:positive}{3} and
\hyperlink{A:density2}{5} hold and $\IK=1$ in Assumption \hyperlink{A:state}{2}.

1. There are $0<c<C<\infty$, such that almost surely, $c\le
\Lambda_{n,a}(\rx)\le C$ for all $n\gg1$, $a\in\cH$ and $\rx$, thus
\[
\mean[\ln\sL_a(\rx)|\eta] =
\lim_{n\toi} \mean[\ln\Lambda_{n,a}(\rx)|\eta].
\]

2. For $\nu$ with $1\le|\nu|\le q$ and $a\in\cH$,
\[
\mean[\ln\sL_a(\rx)|\eta]^{(\nu)}= \mean\bigl[(\ln\sL_a)^{(\nu)}(\rx)|\eta\bigr] =
\lim_{n\toi} \mean\bigl[(\ln\Lambda_{n,a})^{(\nu)}(\rx)|\eta\bigr].
\]

Similar results hold for $\Lambda_{-n,a}$ and $\bar\sL_a$.
\end{theorem}

\section{Binary state HMM with univariate parameters}
\label{sec:example}
In this section, we consider in more detail the case where $\state$ is
a first order binary state Markov chain, with $\state_t = \mathbf{1}\{
H_t \mbox{ is false}\}$. Also, we suppose $\rx\in\Reals$ and
%
%
\begin{equation} \label{eq:degenerate-0}
\prx_0(0) = \prx_1(0)=0,
\end{equation}
that is, at $\rx=0$, false and true nulls are no longer distinguishable
based on the data.

\subsection{Derivatives of likelihood ratio}
We shall focus $t=0$. Analysis for other $t$ can be done likewise.
By (\ref{eq:rho-limit}), the full likelihood ratio (FLR)
$\rho(\rx)$ satisfies
%
%
\begin{equation} \label{eq:FLR-LLR}
\ln\frac{\rho(\rx)}{\tilde\rho(\rx)}
= \sr(\rx) + \bar\sr(\rx)\qquad
\mbox{with }
\sr(\rx) = \ln\frac{\sL_1(\rx)}{\sL_0(\rx)},\
\bar\sr(\rx) =
\ln\frac{\bar\sL_1(\rx)}{\bar\sL_0(\rx)},
\end{equation}
where $\tilde\rho(\rx)$ is the local likelihood
ratio (LLR) for $\state_0$ only based on $X_0$:
\[
\tilde\rho(\rx) =
\frac{\mathsf{Pr}\{\state_0=1\mid X_0\}}{\mathsf{Pr}\{\state_0=0\mid
X_0\}}
=
\frac{P_0(1)\psi_0(\rx,1)}
{P_0(0) \psi_0(\rx,0)},
\]
with $\psi_t(\rx,a)$ being defined in (\ref{eq:psi-t}).

Consider $\sr(\rx)$. The treatment of $\bar\sr(\rx)$ is
similar. By Theorem \ref{thm:L-limit},
%
%
\begin{equation} \label{eq:lambda-n}\quad
\sr(\rx)= \lim_{n\toi} \lambda_n(\rx)\qquad
\mbox{with } \lambda_n(\rx)
=
\ln\frac{
\mean_\stx[\prod_{s=1}^n \psi_s(\rx, \stx_s) \Gv\stx_0=1]
}{
\mean_\stx[\prod_{s=1}^n \psi_s(\rx, \stx_s) \Gv\stx_0=0]
}.
\end{equation}
By (\ref{eq:degenerate-0}), for $t\in\Ints$,
%
%
\begin{equation} \label{eq:degenerate}
\psi_t(0,\stx_t) = f(\tf(Z_t, \prx_{\state_t}(0)),
\prx_{\stx_t}(0))= f(\tf(Z_t, 0), 0)
\end{equation}
is independent of $\stx$, so $\lambda_n(0)=0$, giving $\sr(0)=
0$. Next, define
%
%
\begin{eqnarray} \label{eq:d-D}
d_t(\rx) &=& \ln\psi_t(\rx,1) - \ln\psi_t(\rx,0),\nonumber\\[-8pt]\\[-8pt]
D_{st} &=& P_{st}(1,1) - P_{st}(0,1),\qquad s,
t\in\Ints.\nonumber
\end{eqnarray}

In general, unless $\state$ is stationary, $D_{st}\not=D_{ts}$ for
$s\not=t$. By simple algebra, we have the following identity, which
the next result relies upon
%
%
\begin{equation} \label{eq:Dproduct}
D_{rs} D_{st} = D_{rt},\qquad
D_{ts} D_{sr} = D_{tr},\qquad r\le s \le t.
\end{equation}
\begin{theorem} \label{thm:d-L-1-2}
Let Assumptions \hyperlink{A:noise}{1}--\hyperlink{A:density}{4} hold. Then
%
%
\begin{eqnarray} \label{eq:d-L-1}
\sr'(0)
&=& \sum_{t=1}^\infty D_{0t} d_t'(0), \\
\label{eq:d-L-2}
\sr''(0)
&=& \sum_{t=1}^\infty D_{0t}\{d_t''(0) + [P_{0t}(1,0)-P_{0t}(0,1)]
[d_t'(0)]^2 \}
\nonumber\\[-8pt]\\[-8pt]
&&{}
+ 2\sum_{t=1}^\infty D_{0t} d_t'(0) \sum_{s=1}^{t-1} [P_{0s}(1,0)
- P_{0s}(0,1)] d_s'(0), \nonumber
\end{eqnarray}
where $'$, $'', \ldots,$ denote differentiations with respect to
$\rx$.
\end{theorem}

Simplifications can be made when $\state$ is stationary and
ergodic. In this case, $p_a = P_0(a)\in(0,1)$ and the transition
matrix can be written as
\[
Q = \pmatrix{1\cr1}(p_0, p_1) + r \pmatrix{p_1\cr-p_0} (1, -1),
\]
where $r\in(-1,1)$ is the eigenvalue of $Q$ different from 1. Then
for $t\ge1$,
\[
Q^t = \pmatrix{1\cr1}(p_0, p_1) + r^t
\pmatrix{p_1\cr-p_0}(1, -1)
=
\pmatrix{p_0 + r^t p_1 & p_1 - r^t p_1\cr
p_0 - r^t p_0 & p_1 + r^t p_0},
\]
so that in (\ref{eq:d-L-1}) and (\ref{eq:d-L-2}), $D_{0t} = r^t$ and
$P_{0s}(1,0)-P_{0s}(0,1) = (p_0-p_1)(1-r^s)$.

\subsection{A univariate case} \label{ssec:univariate}
In this section, suppose both $X_t$ and $\prx_{\state_t}(\rx)$ are
univariate. Suppose the following regularity conditions are satisfied:
\begin{enumerate}
\item$\lambda(x,\px) = \ln f(x,\px)\in C^{(2)}$ and $\tf(z,v)$ as a
function in $v$ belongs to $C^{(2)}$, such that for any $\px$, $v$,
and $\nu$ with $|\nu|\le2$,
\[
\mean\bigl[\lambda^{(\nu)}(\tf(Z_t,v), \px)\bigr]=
(\mean[\lambda(\tf(Z_t,v),\px)])^{(\nu)},
\]
where the differentiation is with respect to $v$ and $\px$.
\item$\prx_a(\rx)\in C^{(2)}$ for any $a\in\cH$.
\end{enumerate}
\begin{prop}\label{prop:d-uv}
Let Assumptions \hyperlink{A:noise}{1}--\hyperlink{A:density}{4} hold. Then for
each $t$,
%
%
\begin{eqnarray}\qquad
\label{eq:uv-d-diff-1}
d_t'(0)
&=& [\prx_1'(0)-\prx_0'(0)]\,
\frac{\partial\lambda(x,0)}{\partial\px},
\\
\label{eq:uv-d-diff-2}
d_t''(0)
&=& 2[\prx_1'(0)-\prx_0'(0)] \prx_{\state_t}'(0)\,
\frac{\partial^2 \lambda(x,0)}{\partial x\,\partial\px}\,
\frac{\partial\tf(Z_t,0)}{\partial v}
\nonumber\\[-8pt]\\[-8pt]
&&{}
+[\prx_1'(0)^2-\prx_0'(0)^2]\,
\frac{\partial^2 \lambda(x,0)}{\partial\px^2}
+ [\prx_1''(0)-\prx_0''(0)]\,
\frac{\partial\lambda(x,0)}{\partial\px},\nonumber
\end{eqnarray}
where the partial derivatives of $\lambda$ are evaluated at
$x=\tf(Z_t,0)$.
\end{prop}
\begin{prop} \label{prop:E-d-uv}
Let Assumptions \hyperlink{A:noise}{1}--\hyperlink{A:positive}{3} and
\hyperlink{A:density2}{5} hold and $\IK=1$ in Assumption~\hyperlink{A:state}{2}. Then
%
%
\begin{eqnarray}
\label{eq:E-d-L-1}
\mean[\sr'(0)\mid\state]
&=&0,
\\
\label{eq:E-d-L-2}
\mean[\sr''(0)\mid\state]
&=& \var[d_0'(0)]
\sum_{t=1}^\infty D_{0t}[2\state_t-P_{0t}(1,1)-P_{0t}(0,1)].
\end{eqnarray}
Moreover, for all $t$, $\mean[d_t'(0)]=0$ and
$\var[d_t'(0)]=[\prx_1'(0) - \prx_0'(0)]^2 J(0)$, where $J(\px)$ is
the Fisher information for the parametric family $f(x,\px)$.
\end{prop}

Note that (\ref{eq:E-d-L-2}) implies $\mean[\sr''(0)\mid\state_0]
=(2\state_0-1) \var[d_0'(0)] \sum_{t=1}^\infty D_{0t}^2$,
which is what we got toward the end of Example \ref{ex:gauss}.

We next use the results to get a better view on the structure of
$\rho(\rx)$. Since $\sr(0)=0$, by Taylor's expansion and Theorem
\ref{thm:d-L-1-2},
\begin{eqnarray*}
\sr(\rx)
&=& \sum_{t=1}^\infty D_{0t}
\biggl[d_t'(0) \rx+ \frac{d_t''(0) \rx^2}{2}\biggr]
+
\frac{\rx^2}{2} \sum_{t=1}^\infty D_{0t}
[P_{0t}(1,0)-P_{0t}(0,1)] [d_t'(0)]^2
\\
&&{} + \rx^2\sum_{t=1}^\infty D_{0t} d_t'(0) \sum_{s=1}^{t-1} [P_{0s}(1,0)
- P_{0s}(0,1)] d_s'(0) + o(\rx^2).
\end{eqnarray*}

Since $d_t(0)=0$, then $d_t'(0)\rx+ d_t''(0)\rx^2/2 =
d_t(\rx)+o(\rx^2)$. Under the condition of Proposition
\ref{prop:E-d-uv}, by Propositions \ref{prop:d-uv} and
\ref{prop:E-d-uv}, all $d_t'(0)$ are independent of $\state$, have
mean 0 and the same variance. Similar assertions can be made about
the expansion of $\bar\sr(\rx)$. It follows that
\[
\sr(\rx) + \bar\sr(\rx)
= \sum_{t\not=0} D_{0t} d_t(\rx) + \rx^2 \var[d_0'(0)] K
+ \rx^2 \xi+ o(\rx^2),
\]
where
$K = (1/2) \sum_{t\not=0} D_{0t} [P_{0t}(1,0)-P_{0t}(0,1)]$ and $\xi$
is a random variable independent of $\state$ and has mean 0.
Then by (\ref{eq:FLR-LLR}) and the definition of $d_t(\rx)$ in
(\ref{eq:d-D}),
%
%
\begin{equation} \label{eq:rho-factors}\qquad
\rho(\rx) = \tilde\rho(\rx) \prod_{t\not=0} \biggl[\frac{\psi_t(\rx
,1)}{\psi_t(\rx,0)} \biggr]^{D_{0t}} \times
\exp\bigl\{\rx^2 \bigl(\var[d_0'(0)]K + \xi\bigr) + o(\rx^2)\bigr\}.
\end{equation}

According to (\ref{eq:psi-t}), $\psi_t(\rx,1)/\psi_t(\rx,0)$
is the marginal likelihood ratio of $X_t$ for the isolated test on
$\state_t=1$ vs. $\state_t=0$, which completely ignores the
dependence among the sites. The above expansion shows that
all these likelihood ratios are factored into the FLR, with effects
being adjusted by $D_{0t}$. For example, if $D_{0t}$ is positive
(resp., negative), then a large likelihood ratio at site $t$ increases
(resp., decreases) the FLR for the test on $\state_0$. Also, by
(\ref{eq:Dproduct}), if $s$
has the same sign as $t$ but farther away from 0, then the effect of
the marginal likelihood ratio at site $s$ on the test on $\state_0$ is
determined by $D_{0t}$ and $D_{ts}$. In contrast,
the LLR $\tilde\rho(\rx)$ only takes into account the marginal
likelihood ratio at site 0.

The above expansion is obtained for $\state_0$. Taking into account
explicitly the dependence on site location, the FLRs
for the multiple tests on $\state_s$, $s\in\Ints$, are
%
%
\begin{equation} \label{eq:rho-factors2}\qquad\quad
\rho_s(\rx) = \tilde\rho_s(\rx) \prod_{t\not=s} \biggl[\frac{\psi
_t(\rx,1)}{\psi_t(\rx,0)} \biggr]^{D_{st}} \times
\exp\bigl\{\rx^2 \bigl(\var[d_0'(0)] K_s + \xi_s\bigr) + o(\rx^2)\bigr\},
\end{equation}
where the LLR $\tilde\rho_s(\xi)$ and constants $K_s$ are now
expressed as
\[
\tilde\rho_s(\rx) =
\frac{P_s(1)\psi_s(\rx,1)}{P_s(0)\psi_s(\rx,0)},\qquad
K_s = (1/2) \sum_{t\not=s} D_{st} [P_{st}(1,0)-P_{st}(0,1)]
\]
and $\xi_s$ are centered random variables independent of $\state$.
The conditional likelihoods of $\eta_s$ can then be computed via
$\mathsf{Pr}\{\state_s=0\mid X\} = [1+\rho_s(\rx)]^{-1}$.

\subsection{Examples} \label{ssec:examples}
\begin{example}[(Translation)]\label{ex:shift}
Suppose $\tf$ is defined on $\Reals\times\Reals$ such that
$\tf(z,v)=z+v$ and for $a=0$, $1$, $\prx_a(\rx) = \rx a$. Let each
$Z_t$ have density $h(z) = e^{-V(z)}$. Apparently, Example
\ref{ex:gauss} belongs to this case.

For each $\px\in\Reals$, $\tf(Z_t, \px) = Z_t+\px$ has density
$f(x,\px) = h(x-\px)$. Therefore, $\lambda(x,\px)=\ln f(x,\px) =
-V(x-\px)$. It is easy to check
\begin{eqnarray*}
\theta_a'(0) &=& a,\qquad
\frac{\partial\tf(z,0)}{\partial v} = 1,\qquad
\frac{\partial\lambda(x,\px)}{\partial\px} = V'(x-\px),\\
\frac{\partial^2\lambda(x,\px)}{\partial x\,\partial\px} &=&
- \frac{\partial^2\lambda(x,\px)}{\partial\px^2}=V''(x-\px).
\end{eqnarray*}
Provided necessary conditions are satisfied, by Proposition
\ref{prop:d-uv},
\begin{eqnarray*}
d_t'(0) &=& V'(Z_t),\qquad
d_t''(0) = (2\state_t-1) V''(Z_t),\\
\var[d_t'(0)] &=& \int V'(x)^2 e^{-V(x)} \,dx.
\end{eqnarray*}
Then we can get $\sr'(0)$, $\sr''(0)$ and
$\mean[\sr''(0)\mid\state]$ by Theorem \ref{thm:d-L-1-2} and
(\ref{eq:E-d-L-2}).
\end{example}
\begin{example}[(Scaling)]\label{ex:scaling}
Suppose $\tf$ is defined on $\Reals\times\Reals$ such that
$\tf(z, v) = e^{-v} z$ and for $a=0$, $1$, $\prx_a(\rx) = \rx
a$. Let each $Z_t$ have density $h(z) = e^{-V(z)}$. For
$v\in\Reals$, $\tf(Z_t, v)$ has density $f(x, v)=e^v
h(e^v x)$. Therefore, $\lambda(x,v) = v - V(e^v x)$. By
Proposition \ref{prop:d-uv},
\begin{eqnarray*}
d_t'(0) &=& 1-Z_t V'(Z_t),\qquad
d_t''(0) = (2\state_t-1) Z_t [V'(Z_t) + Z_t V''(Z_t)], \\
\var[d_0'(0)] &=& \int[1-x V'(x)]^2 e^{-V(x)} \,dx.
\end{eqnarray*}
Then we can get $\sr'(0)$, $\sr''(0)$, and
$\mean[\sr''(0)\mid\state]$ by Theorem \ref{thm:d-L-1-2} and
(\ref{eq:E-d-L-2}).
\end{example}
\begin{example}[(\textit{t}-statistics)]\label{ex:t}
Suppose the data observed at each time point $t$ is a random
vector $\xi_t=(\xi_{t,1}, \ldots, \xi_{t,\nu+1})$, such that\vspace*{1pt}
conditional on $\state$, $\xi_t$ are independent of each other, and
for each $t$, $\xi_{t,j}$ are i.i.d. $\sim N(\rx s_t\state_t, s_t^2)$ for
some $s_t = s_t(\state)>0$. Suppose $s_t$ are completely
intractable, that is, there is no information on the values of $s_t$ or
their interrelations. In this case, it is reasonable to use the
$t$-statistics
\[
X_t=\frac{\sqrt{\nu+1}\bar\xi_t}{\sqrt{S_t^2/\nu}}
\]
to test on $\state_t$, where $\bar\xi_t$ is the mean of
$\xi_{t,j}$ and $S_t^2$ the sum of squares of
$\xi_{t,j}-\bar\xi_t$.

By scaling, we
assume without loss of generality that $s_t=1$.
Let $\zeta_t = \sqrt{\nu+1}(\bar\xi_t - \rx\state_t)$. Then, given
$\state$, $\zeta_t\sim N(0,1)$ and $S_t^2\sim\chi_\nu^2$ are
independent of each other. Define $Z_t = (\zeta_t, S_t)$ and, for
$z=(r,s)$ and $a=0,1$, define
\[
\tf(z,v) = \sqrt{\nu}(r+v)/s,\qquad
\prx_a(\rx) = \sqrt{\nu+1} a\rx.
\]
Then $X_t = \sqrt{\nu} (\zeta_t + \sqrt{\nu+1}\state_t\rx)/S_t =
\tf(Z_t, \prx_{\state_t}(\rx))$. Conditional on $\state$, $X_t
\sim t_{\nu, \px}(x)$ with $\px= \prx_{\state_t}(\rx)$, that is,
the noncentral $t$-distribution with $\nu$ degrees of
freedom (df) and noncentrality parameter $\px$. In the notation of
Assumption \hyperlink{A:noise}{1}, $f(x,\px) = t_{\nu,\px}(x)$.

Recall
\begin{eqnarray*}
t_\nu(x) &=& \frac{C_\nu}{(\nu+x^2)^{(\nu+1)/2}}\qquad
\mbox{with }
C_\nu=\frac{\nu^{\nu/2} \Gamma(({\nu+1})/{2})}
{\sqrt{\pi}\Gamma({\nu}/{2})}, \\
t_{\nu,\px}(x) &=& t_\nu(x) e^{-\px^2/2} \Biggl[1+\sum_{k=1}^\infty\frac
{c_k x^k}{(\nu+x^2)^{k/2}} \frac{\px^k}{k!} \Biggr]
\end{eqnarray*}
with $c_k = \frac{\Gamma(({\nu+k+1})/{2})
2^{k/2}}{\Gamma(({\nu+1})/{2})}$.

Therefore,
\[
\lambda(x,\px) = \ln f(x,\px)
= \ln t_\nu(x) - \half\px^2 +
\ln\Biggl[1+\sum_{k=1}^\infty\frac{c_k x^k}{(\nu+x^2)^{k/2}} \frac{\px
^k}{k!} \Biggr].
\]

By $\ln(1+x) = x-\half x^2 + \frac{1}{3} x^3- \cdots,$
\[
\lambda(x,\px)
= \frac{c_1 x}{\sqrt{\nu+x^2}} \px+ \frac{1}{2}
\biggl\{\frac{(c_2 - c_1^2) x^2}{\nu+x^2}-1\biggr\} \px^2
+ \ln t_\nu(x) + O(\px^3).
\]
It follows that
\begin{eqnarray*}
\frac{\partial\lambda(x,0)}{\partial\px} &=&
\frac{c_1 x}{\sqrt{\nu+x^2}},\qquad
\frac{\partial^2\lambda(x,0)}{\partial x\,\partial\px}
= \frac{c_1\nu}{(\nu+x^2)^{3/2}},\\
\frac{\partial^2\lambda(x,0)}{\partial\px^2}
&=&
\frac{(c_2 - c_1^2) x^2}{\nu+x^2}-1.
\end{eqnarray*}

At $\rx=0$, $X_t = \sqrt{\nu}\zeta_t/S_t$. Since $\prx_a'(0) =
\sqrt{\nu+1} a$, (\ref{eq:uv-d-diff-1}) yields
\[
d_t'(0) = \frac{\sqrt{\nu+1} c_1 X_t}{\sqrt{\nu+X_t^2}}
= \frac{\sqrt{2(\nu+1)}\Gamma({\nu}/{2}+1)\zeta_t}
{\Gamma(({\nu+1})/{2})\sqrt{\zeta_t^2 + S_t^2}}.
\]
Next, since $\partial\tf(Z_t,0)/\partial v = \sqrt{\nu}/S_t$, by
(\ref{eq:uv-d-diff-2}),
\[
d_t''(0) = \frac{2c_1(\nu+1)\state_t S_t^2}
{(S_t^2+\zeta_t^2)^{3/2}} + (\nu+1) \biggl[\frac{(c_2-c_1^2) \zeta
_t^2}{S_t^2 + \zeta_t^2} - 1 \biggr].
\]
Then $\sr'(0)$ and $\sr''(0)$ can be calculated by Theorem
\ref{thm:d-L-1-2}.

To apply Proposition \ref{prop:E-d-uv}, we need to check if
Assumption \hyperlink{A:density2}{5} holds. It is not hard to see that
for $g(\rx):=\lambda(\tf(Z_t, \prx_a(\rx)), \prx_b(\rx))$, $g^
{(k)}(\rx)$ is a linear combination of
$S_t^{-j}\,\frac{\partial^j\lambda(x,\px)}{\partial x^j}\,
\frac{\partial^{k-j} \lambda(x,\px)}{\partial\px^{k-j}}$ evaluated
at $x = \tf(Z_t, \prx_a(\rx))$ and $\px=\prx_b(\rx)$. It is also
not hard to see that $\frac{\partial^j\lambda(x,\px)}{\partial x^j}$
and
$\frac{\partial^{k-j} \lambda(x,\px)}{\partial\px^{k-j}}$ are
bounded, so as long as $\mean[S_t^{-jq^2(q+1)/2}]<\infty$ for $j\le
q$, Assumption \hyperlink{A:density2}{5} holds. Since here $q=2$ and
$S_t^2 \sim\chi_\nu^2$, it suffices to have $\nu>12$. Under this
condition,
\[
\var[d_0'(0)] =
\biggl[\frac{\sqrt{2(\nu+1)}\Gamma({\nu}/{2}+1)} {\Gamma(
({\nu+1})/{2})} \biggr]^2
\mean\biggl[\frac{\zeta_t^2}{\zeta_t^2 + S_t^2}\biggr].
\]
Because $S_t^2$ is the sum of squares of $\nu$ i.i.d. $N(0,1)$ random
variables that are independent of $\zeta_t\sim N(0,1)$, by symmetry,
\[
\mean\biggl[\frac{\zeta_t^2}{S_t^2 + \zeta_t^2}\biggr] = \frac{1}{\nu+1}
\quad\implies\quad
\var[d_0'(0)]
=\frac{1}{2} \biggl[\frac{\nu\Gamma({\nu}/{2})}{\Gamma(({\nu
+1})/{2})} \biggr]^2.
\]
Then $\mean[\sr''(0)\mid\state]$ can be calculated by
(\ref{eq:E-d-L-2}).
\end{example}

\section{Technical details} \label{sec:proof}
\subsection{Some inequalities}
For any set $A$, denote by $\#A$ the number of its elements.
\begin{lemma} \label{lemma:ratio-ws}
Let $\cH$ be a finite set and $W_a\ge0$, $V_a\ge0$ for $a\in\cH$
such that $W:=\sum_a W_a>0$ and $V:=\sum_a V_a>0$. Then for any
$x_a$, $a\in\cH$,
\[
\biggl|W^{-1}\sum_a W_a x_a - V^{-1}\sum_a V_a x_a\biggr|
\le{\max_{a,b\in\cH}} |x_a - x_b|
\biggl[1-\frac{\min_a (V_a/W_a)}{\max_a (V_a/W_a)}\biggr].
\]
\end{lemma}
\begin{pf}
Enumerate the elements in $\cH$ in an arbitrary order. Then the
left-hand side equals $|T|/D$, where
\begin{eqnarray*}
T
&
=& \sum_{a,b} (W_a V_b x_a - W_b V_a x_a)
= \sum_{a<b} (W_a V_b - W_b V_a) (x_a-x_b), \\
D
&
=& \sum_{a,b} (W_a V_b + W_b V_a)
\ge\sum_{a<b} (W_a V_b+W_b V_a).
\end{eqnarray*}
Denote $\Delta= \max_{a,b} |x_a-x_b|$. Then
\begin{eqnarray*}
\frac{|T|}{D}
&\le&
\frac{
\Delta\sum_{a<b} |W_a V_b - W_b V_a|
}{
\sum_{a<b} (W_a V_b + W_b V_a)
}
\le
\Delta\max_{a,b} \frac{W_a V_b-W_b V_a}
{W_a V_b+W_b V_a} \\
&
=&
\Delta\biggl[1-\min_{a,b} \frac{2V_a/W_a}{V_a/W_a+V_b/W_b}\biggr]
\le\Delta
\biggl[1-\frac{\min_a (V_a/W_a)}{\max_a (V_a/W_a)}\biggr],
\end{eqnarray*}
completing the proof.
\end{pf}
\begin{lemma} \label{lemma:ratio-ws2}
Let $\cA$ and $\cB$ be finite sets and $W_a$, $V_a$, $x_a>0$
for $a\in\cA\cup\cB$. Then
\[
\biggl|\frac{\sum_{b\in\cB} W_b x_b}{\sum_{a\in\cA} W_a x_a} - \frac
{\sum_{b\in\cB} V_b x_b}{\sum_{a\in\cA} V_a x_a}\biggr|\le\#B \times\biggl(\frac
{\max_{b\in\cB} x_b}{\min_{a\in\cA} x_a}\biggr)
\max_{a\in\cA, b\in\cB}
\biggl|\frac{W_b}{W_a} - \frac{V_b}{V_a}\biggr|.
\]
\end{lemma}
\begin{pf}
The left-hand side equals $|T|/D$, where
\begin{eqnarray*}
T
&=& \sum_{a\in\cA, b\in\cB}
x_a x_b (W_b V_a-W_a V_b) = \sum_{a\in\cA, b\in\cB}
x_a x_b W_a V_a\biggl(\frac{W_b}{W_a} - \frac{V_b}{V_a} \biggr), \\
D
&
=&
\sum_{a, a'\in\cA} x_a x_{a'} W_a V_{a'}
\ge\Bigl(\min_{a\in\cA} x_a\Bigr) \sum_{a\in\cA} W_a V_a
x_a.
\end{eqnarray*}
Then by
\[
|T| \le\#\cB\Bigl(\max_{b\in\cB} x_b\Bigr) \max_{a\in\cA, b\in\cB}
\biggl|\frac{W_b}{W_a} - \frac{V_b}{V_a}\biggr|
\sum_{a\in\cA} W_a V_a x_a,
\]
the lemma follows.
\end{pf}
\begin{lemma} \label{lemma:ws-diff}
Let $\cH$ be a finite set and $q\in\Nats$. For $a\in\cH$, let
$W_a\dvtx\Reals^p\to[0,\infty)$ and $g_a\dvtx\Reals^p\to\Reals$ be $q$
times differentiable. Suppose $W:=\sum_a W_a>0$. Define function
$\bar g = W^{-1}\sum_a W_a g_a$. Enumerate $\cH$ in an arbitrary
order. Then for $\nu$ with $|\nu|=1$,
%
%
\begin{equation} \label{eq:ws-diff1}\quad
\bar g^{(\nu)}
=
W^{-1}\sum_a W_a g_a^{(\nu)}
+ W^{-2} \sum_{a<b} \bigl(W_a^{(\nu)} W_b-W_a W_b^{(\nu)}\bigr) (g_a - g_b),
\end{equation}
and more generally, for $\nu$ with $|\nu|\le q$,
%
%
\begin{equation} \label{eq:ws-diff2}
\bar g^{(\nu)}
=
W^{-1} \sum_a W_a g_a^{(\nu)}+ \sum_{k=2}^{|\nu|+1} \sum_{0\le j <
\nu} W^{-k} U_{k,\nu, j},
\end{equation}
where $U_{k,\nu,j}$ can be written as
\[
U_{k,\nu,j} = \mathop{\sum_{a_1, \ldots, a_{k} \in\cH, a_1 <
a_2}}_{i_1+\cdots+i_k=\nu-j}
c_{k,\nu}(a_1, \ldots, a_{k}, i_1, \ldots, i_{k})
\prod_{s=1}^k W_{a_s}^{(i_s)}\times
\bigl(g_{a_1}^{(j)} - g_{a_2}^{(j)}\bigr),
\]
with $c_{k,\nu}(a_1, \ldots, a_{k}, i_1, \ldots, i_{k})$ being constants.
\end{lemma}
\begin{pf}
If $|\nu|=1$, then
\begin{eqnarray*}
\bar g^{(\nu)}
&=& W^{-1} \sum_a W_a g_a^{(\nu)}+ W^{-1} \sum_a W_a^{(\nu)} g_a -
W^{-2} \sum_a W_a g_a \sum_b W_b^{(\nu)}\\
&=& W^{-1} \sum_a W_a g_a^{(\nu)}
+ W^{-2} \sum_{a\not=b} \bigl(W_a^{(\nu)} W_b-W_a W_b^{(\nu)}\bigr) g_a \\
&=& W^{-1} \sum_a W_a g_a^{(\nu)}
+ W^{-2} \sum_{a<b} \bigl(W_a^{(\nu)} W_b-W_a W_b^{(\nu)}\bigr) (g_a - g_b),
\end{eqnarray*}
showing (\ref{eq:ws-diff1}), and hence (\ref{eq:ws-diff2}) for
$|\nu|=1$. Let $\nu=e+\mu$, where $|e|=1$ and $0\le\mu<\nu$.
Suppose $\bar g{}^{(\mu)}$ has the form (\ref{eq:ws-diff2}). Then
\[
\bar g{}^{(\nu)}= \bigl(\bar g^{(\mu)}\bigr)^{(e)} = \bar f^{(e)}
+
\sum_{k=2}^{|\nu|} \sum_{0\le j<\nu} (W^{-k} U_{k,\mu,j})^{(e)},
\]
where $\bar f = W^{-1} \sum W_a f_a$, with $f_a = g_a^{(\mu)}$. By
(\ref{eq:ws-diff1}),
\begin{eqnarray*}
\bar f^{(e)}
&=& W^{-1} \sum_a W_a f_a^{(e)}
+ W^{-2} \sum_{a<b} \bigl(W_a^{(e)} W_b-W_a W_b^{(e)}\bigr) (f_a - f_b) \\
&=& W^{-1} \sum_a W_a g_a^{(\nu)}
+ W^{-2} \sum_{a<b} \bigl(W_a^{(e)} W_b-W_a W_b^{(e)}\bigr) \bigl(g_a^{(\mu)}-
g_b^{(\mu)}\bigr).
\end{eqnarray*}
On the other hand, for each $k=2,\ldots, |\nu|$ and $0\le j<\nu$,
\[
(W^{-k} U_{k,\mu,j})^{(e)}
= -k W^{-k-1} \sum_{a\in\cH} W_a^{(e)} U_{k,\mu,j}
+ W^{-k} U_{k,\mu,j}^{(e)}.
\]
It is then not hard to see that $\bar g^{(\nu)}$ has the form
(\ref{eq:ws-diff2}). The proof is complete by induction.
\end{pf}

\subsection{Basic facts}
Define matrix-valued functions $L_n(\rx)=(L_{n,a b}(\rx),
a,b \in\cH)$, such that for $n\ge0$,
%
%
\begin{equation} \label{eq:L-end}
L_{\pm n,a b}(\rx) =
\mean_\stx\Biggl[\mathbf{1}\{\stx_{\pm n}=b\} \prod_{s=1}^n\psi_{\pm s}
(\rx,\stx_{\pm s}) \biggm|\stx_0=a \Biggr].
\end{equation}
Then from (\ref{eq:Lambda}),
%
%
\begin{equation} \label{eq:Lambda-L}
\Lambda_{n,a}(\rx) = \frac{
\sum_{b\in\cH} L_{n, a b}(\rx)
}{
\sum_{b\in\cH} L_{n, \imath b}(\rx)
}.
\end{equation}
For ease of notation, when there is no confusion, $\rx$ will be
omitted.
\begin{lemma} \label{lemma:positive}
Let Assumptions \hyperlink{A:noise}{1}--\hyperlink{A:density}{4} hold. Then
for each $n$ and $a$, $b\in\cH$, $L_{n,a b}\in C^{(q)}$, and for
$|n|\ge\IK$, $L_{n,a b}$ is positive and finite.
\end{lemma}
\begin{pf}
By Assumption \hyperlink{A:density}{4}, $\psi_n(\rx, a)\in C^{(q)}$ for each
$n\in\Ints$ and $a\in\cH$, implying $L_{\pm n,a b}\in C^{(q)}$. For
$n\ge\IK$ and $a$, $b\in\cH$, as $P_{0n}(a,b)>0$, there is at least
one $v =(v_1, \ldots, v_{n})$ with $v_n=b$ and $\mathsf{Pr}\{\stx
_1=v_1, \ldots, \stx_n = v_n\mid\stx_0=a\}>0$. For each such $v$ and
$t=1,\ldots,
n$, by Assumption \hyperlink{A:positive}{3}, $\psi_t(\rx,v_t) \in
(0,\infty)$. Therefore, $L_{n,a b}(\rx) \in(0,\infty)$. The proof
for $L_{-n,a b}$ is similar.
\end{pf}

According to the lemma, $\Lambda_{n,a}\in(0,\infty)$ once $|n| \ge
\IK$. Also, by assumptions \hyperlink{A:state}{2} and \hyperlink{A:positive}{3},
$P_0(a)>0$, $\psi_0(\rx,a)>0$. Therefore, $\rho_{m n}(\rx) \in
(0,\infty)$.

The following relation will be repeatedly used:
%
%
\begin{equation} \label{eq:IL0}
L_{n,a b} = \psi_n(\rx,b) \sum_e L_{n-k,a e} I_{n,e b}^{(k)},\qquad
a, b\in\cH, 1\le k < n,
\end{equation}
where
%
%
\begin{equation}\label{eq:IL}
I_{n,e b}^{(k)} = I_{n,e b}^{(k)}(\rx)
=
\mean_\sigma\Biggl[\mathbf{1}\{\stx_n=b\} \prod_{n-k+1}^{n-1} \psi
_t(\rx, \stx_t) \biggm|\stx_{n-k}=e \Biggr].
\end{equation}
%
Similar relation holds when $n<0$.

\subsection{\texorpdfstring{Proof of Theorem
\protect\ref{thm:L-limit}}{Proof of Theorem 2.2}}

\begin{lemma} \label{lemma:rate}
Let Assumptions \hyperlink{A:noise}{1}--\hyperlink{A:positive}{3} hold.

1. Given $a$, $b\in\cH$ and $ \rx$, for $|n|\ge\IK$, $\min_e
\frac{L_{n,be}(\rx)} {L_{n,a e}(\rx)}$ is strictly positive and
increasing in $n$, and $\max_e \frac{L_{n,be}(\rx)} {L_{n,a e}
(\rx)}$ is finite and decreasing in $|n|$.

2. There is an increasing deterministic function $r(\rx_0)\in
(0,1)$, such that given $\rx_0>0$, for almost all realizations of
$Z$ and $\state$,
%
%
\begin{equation}\label{eq:LL}
\Delta_n(\rx)
:=\max_{a,b,c,d} \biggl|\frac{L_{n,b c}(\rx)}{L_{n,a c}(\rx)} - \frac
{L_{n,b d}(\rx)}{L_{n,ad}(\rx)}\biggr|\le C r(\rx_0)^{|n|},\qquad |n|\ge
\kappa, |\rx|\le\rx_0,\hspace*{-30pt}
\end{equation}
where $C=C(\rx_0, Z)$ is a random variable that only depends on
$\rx_0$ and $Z$ and is finite almost surely. Additionally, for
fixed $\rx$, $\Delta_{\pm n}(\rx)$ are decreasing in $n$.
\end{lemma}
\begin{pf}
We only consider $n>0$. The case $n<0$ is similar. Given $a\not=
b\in\cH$, for $n\ge\IK$ and $c\in\cH$, by Lemma
\ref{lemma:positive}, $\frac{L_{n,bc}}{L_{n,ac}}\in(0,\infty)$.
Then by (\ref{eq:IL0}),
%
%
\begin{equation}\label{eq:R-recur}
\frac{L_{n,b c}}{L_{n,a c}} =
\frac{\sum_e L_{n-k,be} I_{n,e c}^{(k)}}{\sum_e L_{n-k,a e}
I_{n,e c}^{(k)}}.
\end{equation}
Letting $k=1$, it is easy to see that
\[
\min_e \frac{L_{n-1,be}}{L_{n-1,a e}}
\le\frac{L_{n,b c}}{L_{n,a c}}
\le\max_e \frac{L_{n-1,be}}{L_{n-1,a e}}\qquad \mbox{all } c\in\cH,
\]
which implies part 1.

Given $1\le k<n$ and $\rx$, for each $a$, $b$, $c$, $d\in\cH$, apply
Lemma \ref{lemma:ratio-ws} to $x_e = \frac{L_{n-k,be}}{L_{n-k,a e}}$,
$W_e = L_{n-k,a e} I_{n,ec}^{(k)}$ and $V_e = L_{n-k,a e}
I_{n,ed}^{(k)}$. Then by (\ref{eq:R-recur}),
\[
\biggl|\frac{L_{n,b c}}{L_{n,a c}}- \frac{L_{n,b d}}{L_{n,ad}}\biggr|\le
\max_{c,d}
\biggl|\frac{L_{n-k,b c}}{L_{n-k,a c}}- \frac{L_{n-k,b d}}{L_{n-k,ad}}
\biggr|\times
\biggl[1-\frac{\min_e I_{n,ed}^{(k)}/ I_{n,e c}^{(k)}} {\max_e
I_{n,ed}^{(k)} / I_{n,e c}^{(k)}} \biggr].
\]
Take maximum over $c$ and $d$ and then over $a$ and $b$. It follows
that
%
%
\begin{equation} \label{eq:d-R}\qquad
\Delta_n(\rx) \le\gamma_n \Delta_{n-k}(\rx)\qquad
\mbox{with }
\gamma_n=\gamma_n(\rx,k)=
1-\frac{\min_{c,d,e} I_{n,ed}^{(k)} / I_{n,e c}^{(k)}}
{\max_{c,d,e} I_{n,ed}^{(k)} / I_{n,e c}^{(k)}}.
\end{equation}

For $z=(z_1, \ldots, z_{\IK-1})\in\cZ^{\IK-1}$, where $\IK$ is as
in Assumption
\hyperlink{A:state}{2}, define
\begin{eqnarray*}
\alpha(z, \rx)
&=& \mathop{\min_{u_t,v_t\in\cH}}_{1\le t\le\IK-1}
\prod_{t=1}^{\IK-1} f(\tf(z_t, \prx_{u_t}(\rx)),
\prx_{v_t}(\rx)),\qquad
\alpha_*(z, \rx_0)
= \inf_{|\rx|\le\rx_0} \alpha(z, \rx), \\
\beta(z, \rx)
&=&\mathop{\max_{u_t,v_t\in\cH}}_{1\le t\le\IK-1}
\prod_{t=1}^{\IK-1} f(\tf(z_t, \prx_{u_t}(\rx)),
\prx_{v_t}(\rx)),\qquad
\beta^*(z, \rx_0)
= \sup_{|\rx|\le\rx_0} \beta(z, \rx).
\end{eqnarray*}
For $n\ge\IK$, let
\begin{eqnarray*}
\zeta_n
&
=& \zeta_n(\rx_0) = \alpha_*(Z_{n-\IK+1}, \ldots, Z_{n-1}, \rx_0),
\\
\xi_n
&
=& \xi_n(\rx_0) = \beta^*(Z_{n-\IK+1}, \ldots, Z_{n-1}, \rx_0).
\end{eqnarray*}

Since $\psi_t(\rx, \stx_t) = f(\tf(Z_t, \prx_{\state_t}(\rx)),
\prx_{\stx_t}(\rx))$, then for $|\rx|\le\rx_0$,
%
%
\begin{eqnarray}
\label{eq:sandwich}
&&\zeta_n \le\prod_{n-\IK+1}^{n-1} \psi_t(\rx,\stx_t) \le\xi_n
\\
\label{eq:I-low}
&&\quad\implies\quad
\zeta_n P_{n-\IK, n}(e,c) \le I_{n,e c}^{(\IK)}(\rx)
\le\xi_n P_{n-\IK,n}(e,c).
\end{eqnarray}

Given $z\in\cZ^{\IK-1}$, by $\#H<\infty$ and Assumption
\hyperlink{A:positive}{3}, $\alpha(z,\rx)$ and $\beta(z,\rx)$ are continuous
in $\rx$ and $0<\alpha(z,\rx)\le\beta(z,\rx)<\infty$, yielding $0<
\alpha_*(z, \rx_0)\le\beta^*(z, \rx_0)<\infty$. As a result,
$\mathsf{Pr}\{0<\zeta_n\le\xi_n <\infty\}=1$. Fix $0<x<y<\infty$,
such that
$p_0:=\mathsf{Pr}\{x\le\zeta_\IK\le\xi_\IK\le y\}>0$. Note that
$x$ and
$y$ can be chosen in such as way that they only depend on $\rx_0$, the
distribution of $Z$, and $\IK$. Because $Z_t$ are i.i.d., from the
definitions of $\zeta_n$ and $\xi_n$, almost surely, there is an
infinite sequence $n_s=n_s(Z,\rx_0)\ge\IK$, $s\ge0$, such that
%
%
\begin{equation} \label{eq:bounds}
x\le\zeta_{n_s}\le\xi_{n_s}\le y
\end{equation}
and furthermore, $n_s$ can be chosen in such a way that
%
%
\begin{equation}\label{eq:N}
n_s\ge n_{s-1}+\IK,\qquad
|\{s\dvtx n_s\le n\}| \ge\frac{p_0 n}{2\IK}\qquad
\mbox{for } n\gg1.
\end{equation}

On the other hand, since $\#\cH>1$, Assumption \hyperlink{A:state}{2} implies
that
%
%
\begin{equation}\label{eq:I-up}
\phi_* \le P_{n-\IK,n}(e,c) \le1-\phi_*\qquad
\mbox{all } c, e\in\cH.
\end{equation}

Combine (\ref{eq:I-low}), (\ref{eq:bounds}) and (\ref{eq:I-up}) to get
\[
0<\phi_*x \le I_{n_s,e c}^{(\IK)}(\rx)\le(1-\phi_*)
y <\infty\qquad
\forall c, e\in\cH
\]
and hence
%
%
\begin{eqnarray}\label{eq:gamma-r}
\gamma_{n_s}&=&1-\frac{\min_{c,d,e} I_{n_s,ed}^{(\IK)}/ I_{n_s,e
c}^{(\IK)}} {\max_{c,d,e} I_{n_s,ed}^{(\IK)}/ I_{n_s,e c}^{(\IK)}}
\le r_0=r_0(\rx_0)\nonumber\\[-8pt]\\[-8pt]
:\!&=&1-\biggl[\frac{\phi_*x}{(1-\phi_*)y}\biggr]^2 < 1.\nonumber
\end{eqnarray}

Now by (\ref{eq:d-R}), $\Delta_{n_s}(\rx)
\le\Delta_{n_s-\IK}(\rx) r_0$. Since $n_{s-1} \le n_s-\IK$ while
(\ref{eq:d-R}) implies that $\Delta_n(\rx)$ is decreasing,
$\Delta_{n_s} (\rx) \le\Delta_{n_{s-1}}(\rx) r_0$ and hence
$\Delta_{n_s} (\rx) \le\Delta_{n_1}(\rx) r_0^{s-1}$ by induction.
For any $n$, if $n_s\le n<n_{s+1}$, then $\Delta_n(\rx) \le\break
\Delta_{n_1}(\rx) r_0^{s-1} \le\Delta_\IK(\rx) r_0^{s-1}$. Combining
(\ref{eq:N}), for $n\gg1$,
\[
\Delta_n(\rx) \le[\Delta_\IK(\rx)/r_0] r(\rx_0)^n\qquad
\mbox{with } r(\rx_0) = r_0^{{p_0}/({2\IK})}.
\]

Note $\Delta_\IK(\rx) \le\max_{a,b,c} \frac{L_{\IK,ac}(\rx)}
{L_{\IK,bc}(\rx)}$. Using (\ref{eq:L-end}) and (\ref{eq:sandwich})
followed by assumption~\hyperlink{A:state}{2}, it is seen that
\[
\max_{a,b,c}
\frac{L_{\IK,a c}(\rx)}{L_{\IK,b c}(\rx)}
\le\frac{\xi_\IK}{\zeta_\IK} \max_{a,b,c}
\frac{P_{0\IK}(b,c)}{P_{0\IK}(a,c)}
\le\frac{(1-\phi_*)\xi_\IK}{\phi_*\zeta_\IK} < \infty.
\]
Therefore, (\ref{eq:LL}) is proved.

To make $r(\rx_0)$ increasing, replace $r(\rx_0)$ with, say,
$[\inf_{c\ge\rx_0} r(c)+1]/2$. From the construction, $r(\rx_0)$
only depends on the distributional properties of $Z$ and $\state$, but
not specific realizations of the processes. Therefore, $r(\rx_0)$ is
deterministic.
\end{pf}
\begin{lemma}\label{lemma:D-Lambda}
Fix $a\in\cH$ and $\rx$.

1. For $a\in\cH$,
\[
0<\inf_{|n|\ge\kappa} \Lambda_{n,a}(\rx) \le\sup_{|n|\ge\kappa}
\Lambda_{n,a}(\rx)<\infty.
\]

2. For $s\ge n\ge\IK$ and $s\le n\le-\IK$,
\[
|\Lambda_{n,a}(\rx) - \Lambda_{s,a}(\rx)| \le
2\Delta_n(\rx)+\Delta_s(\rx).
\]
\end{lemma}
\begin{pf}
From (\ref{eq:Lambda-L}), for $s\ge n\ge\IK$ and $s\le n\le
-\IK$,
\[
\Lambda_{n,a}(\rx),
\frac{L_{s,a e}(\rx)}{L_{s,\imath e}(\rx)}
\in\biggl[\min_e \frac{L_{n,a e}(\rx)}{L_{n,\imath e}(\rx)}, \max_e
\frac{L_{n,a e}(\rx)}{L_{n,\imath e}(\rx)} \biggr].
\]
Together with part 1 of Lemma \ref{lemma:rate}, this yields the first
part of the lemma and also
\[
\biggl|\Lambda_{n,a}(\rx) - \frac{L_{n,a b}(\rx)}{L_{n,\imath b}(\rx)}\biggr|\le
\Delta_n,\qquad
\biggl|\frac{L_{n,a b}(\rx)}{L_{n,\imath b}(\rx)}- \frac{L_{s,a b}(\rx
)}{L_{s,\imath b}(\rx)}\biggr|\le\Delta_n,
\]
where $b\in\cH$ is arbitrary. Then by
\begin{eqnarray*}
|\Lambda_{n,a}(\rx)- \Lambda_{s,a}(\rx)|&\le&
\biggl|\Lambda_{n,a}(\rx) - \frac{L_{n,a b}(\rx)}{L_{n,\imath b}(\rx)}\biggr|+
\biggl|\Lambda_{s,a}(\rx)- \frac{L_{s,a b}(\rx)}{L_{s,\imath b}(\rx)}\biggr|\\
&&{}+
\biggl|\frac{L_{n,a b}(\rx)}{L_{n,\imath b}(\rx)}- \frac{L_{s,a b}(\rx
)}{L_{s,\imath b}(\rx)} \biggr|,
\end{eqnarray*}
the second part of the lemma follows.
\end{pf}

\begin{pf*}{Proof of Theorem \ref{thm:L-limit}}
From Lemmas \ref{lemma:rate}--\ref{lemma:D-Lambda}, it is seen
that given $\rx_0>0$, almost surely, as $n\toi$,
$\Lambda_{n,a}(\rx)\to\sL_a(\rx)$ and
$\Lambda_{-n,a}(\rx)\to\bar\sL_a(\rx)$ uniformly for
$|\rx|\le\rx_0$, at rate $o(r(\rx_0)^n)$. Since $\Lambda_{\pm
n,a}(\rx)$ are continuous, the uniform convergence implies that
$\sL_a(\rx)$ and $\bar\sL_a(\rx)$ are continuous. Also, the
lemmas imply that $\sL_a(\rx)$ and $\bar\sL_a(\rx)$ are strictly
positive. By monotonicity argument, almost surely, the
convergence holds simultaneously for all $\rx_0>0$.
\end{pf*}

\subsection{\texorpdfstring{Proof of Theorem
\protect\ref{thm:d-L-limit}}{Proof of Theorem 2.3}}
For $t\not=0$, $n\ge1$ and $\rx_0>0$, define
%
%
\begin{eqnarray}\label{eq:V-n-def}
V_{\pm n}(\rx_0) &=&
n \max_{1\le t\le n} D_{\pm t}(\rx_0)\hspace*{60pt}
\nonumber\\[-8pt]\\[-8pt]
&&\eqntext{\mbox{with }
\displaystyle D_t(\rx_0) = \max_{|\nu|\le q}\max_{a\in\cH}
\sup_{|\rx|\le\rx_0}
\biggl|\frac{\psi^{(\nu)}_t(\rx, a)}{\psi_t(\rx,a)} \biggr|,}
\end{eqnarray}
where $\psi^{(\nu)}_t$ is a derivative with respect to $\rx$. Note
$D_t(\rx_0)\ge1$ since the maximization in its definition takes into
account $\nu=0$.
\begin{lemma} \label{lemma:bound}
The following statements are true.

1. For $\rx_0>0$ and $n\ge1$,
%
%
\begin{equation}\label{eq:V-n-bound}
V_n(\rx_0) \le n \max_{|t|\le n} [q+M_q(Z_t,\rx_0)]^q.
\end{equation}

2. If Assumptions \hyperlink{A:noise}{1}--\hyperlink{A:density}{4} hold, then
$\mathsf{Pr}\{\lim_n \beta^{-n} V_n(\rx_0)=0, \forall\beta>1,
\rx_0>0\} = 1$.
\end{lemma}
\begin{pf}
To show part 1, it suffices to show that for all $\nu$ with
$|\nu|=l\le q$, and all $\rx_0>0$ and $t\not=0$,
%
%
\begin{equation}\label{eq:dq}
d_{\nu,t}(\rx_0) := \max_{a\in\cH} \sup_{|\rx|\le\rx_0}
\biggl|\frac{\psi^{(\nu)}_t(\rx, a)}{\psi_t(\rx,a)}\biggr|\le[l+M_l(Z_t,\rx_0)]^l.
\end{equation}

It is easily seen that (\ref{eq:dq}) holds for $l=0$. Suppose
(\ref{eq:dq}) holds if $|\nu|\le l$. Let $|\nu|=l+1$. Without loss
of generality, let $\nu=e+\mu$, where $e=(1,0,\ldots,0)$ and
$\mu=(\mu_1, \ldots, \mu_p)\ge0$. Let $\ell_{z,a b}(\rx) = \ln
f(\tf(z, \prx_a(\rx)), \prx_b(\rx))$ as in Assumption
\hyperlink{A:density}{4}. Then by $\psi_t^{(e)}(\rx,a) = \psi_t(\rx,a)
\ell^{(e)}_{Z_t,\state_ta}(\rx)$,
\[
\psi_t^{(\nu)}(\rx,a)
=
\bigl[\psi_t(\rx,a) \ell^{(e)}_{Z_t,\state_t a}(\rx) \bigr]^{(\mu)}
=
\sum_{i\le\mu} \pmatrix{\mu\cr i} \psi_t^{(i)}(\rx,a)
\ell^{(\nu-i)}_{Z_t,\state_t a}(\rx),
\]
where ${\mu\choose i} = {\mu_1\choose i_1} \cdots
{\mu_p\choose i_p}$.

For $i\le\mu$, $|\ell^{(\nu-i)}_{Z_t,\state_t a}(\rx)| \le
M_l(Z_t,\rx_0)$. Then, as $|\mu|=l$, by induction the hypothesis,
\begin{eqnarray*}
\max_{a\in\cH}\sup_{|\rx|\le\rx_0} \biggl|\frac{\psi_t^{(\nu)}(\rx
, a)}{\psi_t(\rx,a)}\biggr|
&\le&
M_l(Z_t,\rx_0)\sum_{i\le\mu} \pmatrix{\mu\cr i}
[|i|+M_l(Z_t,\rx_0)]^{|i|} \\
&\le&
M_l(Z_t,\rx_0)\sum_{i\le\mu} \pmatrix{\mu\cr i}
[l+M_l(Z_t,\rx_0)]^{|i|}\\
&=&
M_l(Z_t,\rx_0) [|\nu|+M_l(Z_t,\rx)]^l,
\end{eqnarray*}
which implies (\ref{eq:dq}). By induction, (\ref{eq:dq}) holds for
all $|\nu|\le q$.

Because $V_n(\rx_0)$ is increasing in $\rx_0$, to show part 2, it
suffices to show that for each fixed $\rx_0>0$ and
$\beta>1$, $\lim_n \beta^{-n} V_n(\rx_0)=0$ almost surely. Fix an
arbitrary $c\in(1,\beta)$, such that $c^q<\beta$. By part 1 and
Assumption \hyperlink{A:density}{4}, for some $p=p(\rx_0)>2$,
\begin{eqnarray*}
\mathsf{Pr}\{V_n(\rx_0)\ge n c^{qn}\}
&\leq& \mathsf{Pr}\Bigl\{\max_{|t|\le n} M_q(Z_t,\rx_0)\ge c^n
-q\Bigr\}\\
&\le&2 n \operatorname{\mathsf{Pr}}\{M_q(Z_0,\rx_0)\ge c^n-q\}
= o(n^{-p+1}).
\end{eqnarray*}
Then part 2 follows from the Borel--Cantelli lemma and $n c^{qn} =
o(\beta^n)$.
\end{pf}
\begin{lemma} \label{lemma:lambda-d}
Let Assumptions \hyperlink{A:noise}{1}--\hyperlink{A:density}{4} hold. Fix $a$,
$b$, $c\in\cH$ and $k\ge1$. Let
\[
W_n(\rx) = L_{n-k,a b}(\rx) I_{n,b c}^{(k)}(\rx),\qquad n\ge k,
\]
where $I_{n,bc}^{(k)}$ is defined in (\ref{eq:IL}). Given $\nu>0$
with $|\nu|\le q$ and $\rx_0>0$, for $n\ge0$,
\[
\sup_{|\rx|\le\rx_0}
\frac{|L_{n,a b}^{(\nu)}(\rx)|}{L_{n,a b}(\rx)} \le
[V_n(\rx_0)]^{|\nu|},
\]
with $V_n(\rx_0):=0$ if $n=0$, while for $n\ge k$,
\[
\sup_{|\rx|\le\rx_0} \frac{|W^{(\nu)}_n(\rx)|}{W_n(\rx)}
\le[V_{n-1}(\rx_0)]^{|\nu|}.
\]
\end{lemma}
\begin{pf}
For $\nu=(\nu_1, \ldots, \nu_p)$ with $1\le|\nu|\le q$, it is not
hard to
get
\[
L^{(\nu)}_{n,a b}(\rx)
=
\mean_\sigma\Biggl[\mathbf{1}\{\stx_n=b\} \sum_{l_1+\cdots+l_n=\nu}
\prod_{t=1}^n \psi_t^{(l_t)}(\rx, \stx_t) \biggm|\stx_0=a
\Biggr].\vadjust{\goodbreak}
\]

For any sequence $l_1, \ldots, l_{n}$ in the sum, at most $|\nu|$ of them
are nonzero. For each $l_t>0$, $|\psi_t^{(l_t)}(\rx, \stx_t)| \le
D_t(\rx_0) \psi_t(\rx, \stx_t)$ for $|\rx|\le\rx_0$. As a result,
\[
\prod_{t=1}^n \bigl|\psi_t^{(l_t)}(\rx, \stx_t)\bigr|
\le\Bigl[\max_{1\le t\le n} D_t(\rx_0)\Bigr]^{|\nu|}
\prod_{t=1}^n \psi_t(\rx,\stx_t).
\]
On the other hand, there are $n^{\nu_1}\cdots n^{\nu_p} = n^{|\nu|}$
such sequences. Then
\begin{eqnarray*}
\bigl|L^{(\nu)}_{n,a b}(\rx)\bigr|
&\le&
\Bigl[n\max_{1\le t\le n} D_t(\rx_0)\Bigr]^{|\nu|}
\mean_\sigma\Biggl[\mathbf{1}\{\stx_n=b\} \prod_{t=1}^n\psi_t(\rx,
\stx_t) \biggm|\stx_0=a \Biggr] \\
&=&
\Bigl[n\max_{1\le t\le n} D_t(\rx_0)\Bigr]^{|\nu|} L_{n,a b}(\rx).
\end{eqnarray*}
This completes the proof of the first inequality. To show the
second inequality, first,
\[
W_n^{(\nu)}(\rx) = \sum_{i\le\nu} \pmatrix{\nu\cr i}
L_{n-k,a b}^{(i)}(\rx) \bigl[I_{n,b c}^{(k)}\bigr]^{(\nu-i)}(\rx).
\]

Using the definition of $I_{n,bc}^{(k)}$ and following the treatment
for $L_{n,a b}^{(\nu)}(\rx)$,
\[
\bigl|\bigl[I_{n,b c}^{(k)}\bigr]^{(\nu-i)}(\rx) \bigr|\le
(k-1)^{|\nu|-|i|}
\Bigl[\max_{n-k+1\le t\le n-1} D_t(\rx_0)\Bigr]^{|\nu|-|i|} I_{n,b
c}^{(k)}(\rx).
\]
Combining the bound with the one for $L_{n-k,a b}^{(i)}(\rx)$,
\begin{eqnarray*}
\bigl|W_n^{(\nu)}(\rx)\bigr|
&
\le&
\Bigl[\max_{1\le t\le n-1} D_t(\rx_0)\Bigr]^{|\nu|}\\
&&{}\times\sum_{i\le\nu} \pmatrix{\nu\cr i}(n-k)^{|i|} (k-1)^{|\nu|-|i|}
L_{n-k,a b}(\rx) I_{n,b c}^{(k)}(\rx) \\
&
\leq&
[V_{n-1}(\rx_0)]^{|\nu|} W_n(\rx).
\end{eqnarray*}
This finishes the proof.
\end{pf}
\begin{lemma} \label{lemma:diff-rate}
Let Assumptions \hyperlink{A:noise}{1}--\hyperlink{A:density}{4} hold. Define, for
$\nu$ with $|\nu|=1,\ldots, q$,
%
%
\begin{equation}\label{eq:LL-d}
\Delta_{n,\nu}(\rx)
:=\max_{a,b,c,d} \biggl|\biggl(\frac{L_{n,b c}}{L_{n,a c}}\biggr)^{(\nu)} (\rx) -
\biggl(\frac{L_{n,b d}}{L_{n,ad}}\biggr)^{(\nu)} (\rx) \biggr|.
\end{equation}
Then for each $\nu$, there is an increasing deterministic function
$0\le r_\nu(\rx_0)<1$ in $\rx_0>0$, such that almost surely, as
$n\toi$,
\[
\sup_{|\rx|\le\rx_0} \Delta_{n,\nu}(\rx) = o\bigl(r_\nu(\rx_0)^{|n|}\bigr)
\qquad\mbox{all } \rx_0>0.
\]
\end{lemma}
\begin{pf}
We only consider $n>0$. The case $n<0$ can be handled similarly.
Given $k$, define $I_{n,ec}^{(k)}(\rx)$ as in (\ref{eq:IL}). Given
$a\not=b\in\cH$, write $W_{n,ec}(\rx) = L_{n-k,a e}(\rx
)I_{n,ec}^{(k)}(\rx)$, $W_{n,c}(\rx) = \sum_e W_{n,ec}(\rx)$. Then by
(\ref{eq:R-recur}), for $n\ge\IK$,
\[
\frac{L_{n,b c}}{L_{n,a c}}
= W_{n,c}^{-1} \sum_e W_{n,e c} \frac{L_{n-k,be}}{L_{n-k,a e}}.
\]

Fix $l=1,\ldots, q$. By Lemma \ref{lemma:ws-diff}, for $\nu\not=0$
with $|\nu|=l$,
%
%
\begin{equation}\label{eq:d-lambda}
\biggl(\frac{L_{n,b c}}{L_{n,a c}}\biggr)^{(\nu)}= W_{n,c}^{-1}
\sum_e W_{n,e c} \biggl(\frac{L_{n-k,be}}{L_{n-k,a e}}\biggr)^{(\nu)}
+R_{n,\nu,c},
\end{equation}
where
\[
R_{n,\nu,c} =
\mbox{a linear combination of }
\Biggl[\prod_{s=1}^m \frac{W_{n,e_s c}^{(i_s)}}{W_{n,c}}\Biggr]
\biggl[\biggl(\frac{L_{n-k,be_1}}{L_{n-k,a e_1}}\biggr)^{(j)} -\biggl(\frac
{L_{n-k,be_2}}{L_{n-k,a e_2}}\biggr)^{(j)} \biggr]
\]
across $m=2, \ldots, |\nu|+1$, $i_1, \ldots, i_{m}\ge0$, $0\le
j<\nu$ with
$i_1+\cdots+i_m+j=\nu$, and $e_1, \ldots, e_{m}\in\cH$ with
$e_1<e_2$. Then,
by the same argument that leads to (\ref{eq:d-R}),
%
%
\begin{equation} \label{eq:d-R-n}
\Delta_{n,\nu}(\rx) \le\gamma_n \Delta_{n-k,\nu}(\rx) +
{2 \max_c} |R_{n,\nu,c}(\rx)|,
\end{equation}
where $\gamma_n$ is given in (\ref{eq:d-R}).

The rest of the proof is by induction on $l$. First, let
$|\nu|=1$. By Lemma \ref{lemma:ws-diff},
\[
R_{n,\nu,c}=
W_{n,c}^{-2}\sum_{e_1<e_2}
\bigl(W_{n,e_1c}^{(\nu)} W_{n,e_2c} - W_{n,e_1c} W_{n,e_2c}^{(\nu)}\bigr)
\biggl(\frac{L_{n-k,be_1}}{L_{n-k,a e_1}} -\frac{L_{n-k,be_2}}{L_{n-k,a
e_2}} \biggr).
\]

Fix $\rx_0>0$. By Lemma \ref{lemma:lambda-d}, for $|\rx|\le\rx_0$,
$|W_{n,e c}^{(\nu)}(\rx)| \le V_{n-1}(\rx_0) W_{n,e c}(\rx)$. Then
%
%
\begin{eqnarray}\label{eq:R-D}
&&|R_{n,\nu,c}(\rx)|
\le W_{n,c}^{-2} \sum_{e_1<e_2}
2 V_{n-1}(\rx_0) W_{n, e_1 c} W_{n, e_2 c} \Delta_{n-k}(\rx)
\nonumber\\[-8pt]\\[-8pt]
&&\quad\implies\quad
{\max_c}
|R_{n,\nu,c}(\rx)| \le V_{n-1}(\rx_0)\Delta_{n-k}(\rx).\nonumber
\end{eqnarray}
By Lemma \ref{lemma:rate}, there is increasing deterministic
$r=r(\rx_0)\in(0,1)$, such that $\sup_{|\rx|\le\rx_0} \Delta
_n(\rx)
\le r^n$ for $n\gg1$. Fix $\beta\in(1,1/r)$. Then by
(\ref{eq:d-R-n}), (\ref{eq:R-D}) and part~2 of Lemma
\ref{lemma:bound}, almost surely, for $n\gg1$ and $|\rx| \le\rx_0$,
%
%
\begin{equation}\label{eq:R-D2}
\Delta_{n,\nu}(\rx)
\le
\gamma_n\Delta_{n-k,\nu}(\rx) + \beta^n r^{n-k}
\le\Delta_{n-k,\nu}(\rx) + \beta^n r^{n-k}.
\end{equation}
Let $k=1$ to get $\Delta_{n,\nu}(\rx)\le\Delta_{n-1,\nu}(\rx) +
\beta^n r^{n-1}$. So by induction, for $s\le n$,
%
%
\begin{equation}\label{eq:D-recur}
\Delta_{n,\nu}(\rx) \le\Delta_{s,\nu}(\rx) + \beta
\sum_{t=s}^{n-1}(\beta r)^t \le\Delta_{s,\nu}(\rx) +
\frac{\beta}{1-\beta r} (\beta r)^s.
\end{equation}

Next let $k=\IK$. By the same argument that leads to
(\ref{eq:gamma-r}), $r$ can be chosen in such a way that there is a
sequence $n_s=n_s(Z,\rx_0)$ that satisfy (\ref{eq:N}) and
$\gamma_{n_s}\le r$. By the first inequality in (\ref{eq:R-D2}),
for $s\gg1$,
\[
\Delta_{n_s,\nu}(\rx)
\le
r\Delta_{n_s-\IK,\nu}(\rx) + \beta^{n_s} r^{n_s-\IK}
\le
r\Delta_{n_s-\IK,\nu}(\rx) + \beta^\IK(\beta r)^{n_{s-1}}.
\]
Let $n=n_s-\IK$ and $s=n_{s-1}$ in (\ref{eq:D-recur}) and combine it
with the last equality to get
\[
\Delta_{n_s,\nu}(\rx) \le r\Delta_{n_{s-1},\nu}(\rx) + c(\beta
r)^{n_{s-1}},
\]
where $c=\beta^\IK+ \beta/(1-r\beta)$. Then by induction and the
fact that $n_s\ge\IK s$,
\begin{eqnarray*}
\Delta_{n_s,\nu}(\rx)
&
\le& r^{s-1} \Delta_{n_1,\nu}(\rx) + c
\sum_{t=1}^{s-1} r^{s-t-1} (\beta r)^{n_t} \\
&
\le&
r^{s-1} \Delta_{n_1,\nu}(\rx) + c
\sum_{t=1}^{s-1} r^{s-t-1} (\beta r)^t\\
&\le&
r^{s-1} \Delta_{n_1,\nu}(\rx) + c s (\beta r)^{s-1}.
\end{eqnarray*}
Now for any $n_s\le n< n_{s+1}$, by (\ref{eq:D-recur}) and the above
inequality,
\[
\Delta_{n,\nu}(\rx) \le r^{s-1} \Delta_{n_1,\nu} +
\biggl(\frac{\beta}{1-r\beta} + c s\biggr)(\beta r)^{s-1}.
\]
Since for $s\gg1$, $s+1\ge\frac{p_0}{2\IK} n_{s+1} \ge\frac{p_0}
{2\IK} n$, it can be seen that $\Delta_{n,\nu}(\rx) = O(c^n)$, with
$c= (\beta r)^{{p_0}/({2\IK})}<1$. Since $\beta\in(1,1/r)$ is
arbitrary, it follows that for a given $\rx_0$ and any $r_1 > r_* :=
r^{{p_0}/({2\IK})}$, say $r_1=r_1(\rx_0) = (1+r_*)/2$,
$\sup_{|\rx|\le\rx_0} \Delta_{n,\nu}(\rx)=o(r_1^n)$ almost surely.
By monotonicity, it can be seen that almost surely, the
exponentially fast convergence holds simultaneously for all
$\rx_0$.

Now let $|\nu|>1$. To bound $R_{n,\nu,c}(\rx)$, for $s=2, \ldots,
|\nu|+1$, and $p$-tuples of nonnegative integers, $i_1, \ldots,
i_{s}$, $j$,
$i_1+\cdots+i_s =\nu-j<\nu$, and $e_1, \ldots, e_{s}\in\cH$, by Lemma
\ref{lemma:lambda-d}, for $|\rx|\le\rx_0$,
\[
\bigl|W_{n,e_1c}^{(i_1)}\cdots W_{n,e_s c}^{(i_s)}\bigr|\le\prod_{k=1}^s
[V_{n-1}(\rx_0)]^{|i_k|} W_{n,e_k c}
\le W_{n,c}^s [V_{n-1}(\rx_0)]^{|\nu|}
\]
so in place of (\ref{eq:R-D}), we have
%
%
\begin{equation} \label{eq:R-D-n}
{\max_c} |R_{n,\nu,c}(\rx)| \le
C_\nu[V_{n-1}(\rx_0)]^{|\nu|}
\times\sum_{j<\nu} \Delta_{n-k,j}(\rx),
\end{equation}
where $\Delta_{n-k,0}(\rx):=\Delta_{n-k}(\rx)$ and $C_\nu>0$ is some
constant only depending on $\nu$.

Suppose it has been shown that for each $j<\nu$, there is $r_j =
r_j(\rx_0)\in(0,1)$, such that $\sup_{|\rx|\le\rx_0}\Delta_{n,j}
(\rx) = o(r_j^n)$. Then using\vspace*{1pt} (\ref{eq:d-R-n}) and (\ref{eq:R-D-n})
and following the argument for $\Delta_{n,j}(\rx)$ with $|j|=1$,
$\sup_{|\rx|\le\rx_0} \Delta_{n,\nu}(\rx) = o(r_\nu^n)$ for some
$r_\nu= r_\nu(\rx_0)\in(0,1)$. By induction, the exponential rate
of convergence holds for all $\nu$ with $|\nu|\le q$. Again, from
the construction, $r_\nu$ only depends on the distributional
properties of $Z$ and $\state$ and hence is deterministic.
\end{pf}

Set $k=1$ in (\ref{eq:d-lambda}). For $n\ge\IK$ and $a$, $b$,
$c\in\cH$,
\[
\min_e\biggl(\frac{L_{n-1,be}}{L_{n-1,a e}}\biggr)^{(\nu)}-|R_{n,\nu,c}|
\le
\biggl(\frac{L_{n,b c}}{L_{n,a c}}\biggr)^{(\nu)}
\le
\max_e\biggl(\frac{L_{n-1,be}}{L_{n-1,a e}}\biggr)^{(\nu)}+|R_{n,\nu,c}|,
\]
giving
%
%
\begin{equation}\label{eq:L-ratio-recurse}
\biggl|\biggl(\frac{L_{n,b c}}{L_{n,a c}}\biggr)^{(\nu)}(\rx) - \biggl(\frac{L_{n-1,b
c}}{L_{n-1,a c}}\biggr)^{(\nu)}(\rx)\biggr|\le\Delta_{n-1,\nu}(\rx)+2 |R_{n,\nu
,c}(\rx)|.
\end{equation}
\begin{cor} \label{cor:diff-rate}
Let assumptions \hyperlink{A:noise}{1}--\hyperlink{A:density}{4} hold. Then almost
surely, as $s\ge n\toi$,
\begin{eqnarray*}
{\max_{a\in\cH}\sup_{|\rx|\le\rx_0}} |R_{n,\nu,c}(\rx)| &=&
o(r_\nu^n(\rx_0)) \\
\max_{a,b,c\in\cH} \sup_{|\rx|\le\rx_0} \biggl|\biggl(\frac{L_{n,b
c}}{L_{n,a c}}\biggr)^{(\nu)} (\rx) - \biggl(\frac{L_{s,b c}}{L_{s,a c}}\biggr)^{(\nu
)} (\rx)\biggr|&=& o(r_\nu^n(\rx_0)),
\end{eqnarray*}
for all $\rx_0>0$ and $\nu$ with $1\le|\nu|\le, q$, and likewise
for $\bar L_{n,a b}$, where $r_\nu(\rx_0)$ are given in Lemma
\ref{lemma:diff-rate}.
\end{cor}
\begin{pf}
The first inequality is already shown in the proof of Lemma
\ref{lemma:diff-rate}. The second one follows from summing the
inequality in (\ref{eq:L-ratio-recurse}) over $n+1, \ldots, s$ and
applying the first inequality and Lemma \ref{lemma:diff-rate}.
\end{pf}

\begin{pf*}{Proof of Theorem \ref{thm:d-L-limit}}
Let $r_\nu(\rx_0)$ be as in Lemma \ref{lemma:diff-rate}. For
$e\in\cH$, denote
\[
\omega_{n,e} = L_{n,\imath e},\qquad
\omega_n = \sum_e \omega_{n,e}.
\]
Then for $a\in\cH$, $\Lambda_{n,a} = \omega_n^{-1}\sum_e
\omega_{n,e} (\frac{L_{n,a e}}{L_{n,\imath e}})$ and similar to
(\ref{eq:d-lambda}),
%
%
\begin{equation} \label{eq:d-Lambda}
\Lambda_{n,a}^{(\nu)}
= \omega_n^{-1}
\sum_e\omega_{n,e} \biggl(\frac{L_{n,a e}}{L_{n,\imath e}}\biggr)^{(\nu)}
+ T_{n,\nu},
\end{equation}
where $T_{n,\nu}$ is a linear combination of
\[
\omega_n^{-m} \omega_{n,e_1}^{(i_1)}\cdots\omega_{n,e_m} ^{(i_m)}
\biggl[\biggl(\frac{L_{n-k,a e_1}}{L_{n-k,\imath e_1}} \biggr)^{(j)} - \biggl(\frac{L_{n-k,a
e_2}}{L_{n-k,\imath e_2}}\biggr)^{(j)} \biggr]
\]
across $m=2, \ldots, |\nu|+1$, $0\le j<\nu$, $i_1, \ldots, i_{m}\ge
0$ with
$i_1+\cdots+i_m+j=\nu$, and $e_1, \ldots, e_{m}\in\cH$ with
$e_1<e_2$. Fix
any $b\in\cH$. From the above formulas,
%
%
\begin{equation} \label{eq:d-Lambda-2a}
\biggl|\Lambda_{n,a}^{(\nu)}-\biggl(\frac{L_{n,a b}}{L_{n,\imath
b}}\biggr)^{(\nu)}\biggr|
\le\Delta_{n,\nu}+|T_{n,\nu}|.
\end{equation}
Following the treatment of $R_{n,\nu,c}$ in (\ref{eq:R-D-n}), except
that we have to use the first inequality in Lemma
\ref{lemma:lambda-d}, it can be seen that
%
%
\begin{equation} \label{eq:T-D-n}
|T_{n,\nu}(\rx)|\le
C_\nu[V_n(\rx_0)]^{|\nu|}
\times\sum_{j<\nu} \Delta_{n-k,j}(\rx),\qquad
|\rx|\le\rx_0,
\end{equation}
yielding $\max_{|\rx|\le\rx_0} |T_{n,\nu}(\rx)|=o(r_\nu^n)$. Now for
$s\not=n$, by (\ref{eq:d-Lambda-2a}), it is not hard to get
%
%
\begin{eqnarray}\label{eq:d-Lambda-2}
\bigl|\Lambda_{s,a}^{(\nu)}-\Lambda_{n,a}^{(\nu)}\bigr|&\le&
\Delta_{s,\nu}+|T_{s,\nu}|+\Delta_{n,\nu}+|T_{n,\nu}|\nonumber\\[-8pt]\\[-8pt]
&&{}
+\biggl|\biggl(\frac{L_{s,a b}}{L_{s,\imath b}}\biggr)^{(\nu)} -\biggl(\frac{L_{n,a
b}}{L_{n,\imath b}}\biggr)^{(\nu)} \biggr|.\nonumber
\end{eqnarray}
Then by Lemma \ref{lemma:diff-rate} and Corollary \ref{cor:diff-rate},
\[
\sup_{|\rx|\le\rx_0}
\bigl|\Lambda_{s,a}^{(\nu)}(\rx) -\Lambda_{n,a}^{(\nu)}(\rx)\bigr|
= o(r_\nu^{s\wedge n}(\rx_0)),\qquad a\in\cH.
\]
Since $\#\cH<\infty$, almost surely, the rate holds simultaneously
for all $a\in\cH$.
\end{pf*}

\subsection{\texorpdfstring{Proof of Theorem \protect\ref{thm:k=1}}{Proof of Theorem 2.4}}
Since the parameter $\IK$ in Assumption \hyperlink{A:state}{2}
equals~1, $P_{n-1,n}(a,b) \in[\phi_*,1-\phi_*]$ for $a$,
$b\in\cH$ and $n\in\Ints$, with $0<\phi_*<1$ as in Assumption
\hyperlink{A:state}{2}. Consequently,
%
%
\begin{equation} \label{eq:k=1-gamma}
\gamma= 1-\inf_n
\frac{
\min_{c,d,e} (({P_{n-1,n}(e,d)})/({P_{n-1,n}(e,c)}))
}
{
\max_{c,d,e} (({P_{n-1,n}(e,d)})/({P_{n-1,n}(e,c)}))
} \in
\biggl[0,1-\frac{\phi_*}{1-\phi_*}\biggr].\hspace*{-35pt}
\end{equation}

For $a$, $e\in\cH$, by (\ref{eq:L-end}), $L_{1,a
e}(\rx) = P_{01}(a,e) \psi_1(\rx,e)$, giving
%
%
\begin{equation} \label{eq:k=1-L-const}
\frac{L_{1,b e}(\rx)}{L_{1,a e}(\rx)}
\equiv
\frac{P_{01}(b,e)}{P_{01}(a,e)}
\le\frac{1-\phi_*}{\phi_*}\qquad
\forall\rx.
\end{equation}
Then by Lemma \ref{lemma:rate},
%
%
\begin{equation} \label{eq:Lambda-ul}
\frac{\phi^*}{1-\phi^*} \le\Lambda_{n,a}(\rx)\le
\frac{1-\phi^*}{\phi^*}.
\end{equation}
This together with dominated convergence shows part 1 of Theorem \ref{thm:k=1}. To prove part 2, we
need several lemmas.
\begin{lemma} \label{lemma:Lambda-step}
Fix $\rx_0>0$. Let $\gamma$ and $\phi_*$ be as in
(\ref{eq:k=1-gamma}) and $\alpha=\phi_*^{-1}-1$. There is a
constant $C>0$, such that if $1\le|\nu|=l\le q$, $|\rx|\le\rx_0$
and $n\ge1$, then
%
%
\begin{eqnarray} \label{eq:k=1-step}
&&\bigl|\Lambda_{n,a}^{(\nu)}(\rx) -
\Lambda_{n-1,a}^{(\nu)}(\rx)\bigr|\nonumber\\[-8pt]\\[-8pt]
&&\qquad
\le C \alpha\gamma^{(n-l-1)\vee1} n^{l(l+2)}
\sum_{t=1}^n [q+M_q(Z_t,\rx_0)]^{q l(l+1)/2}.\nonumber
\end{eqnarray}
\end{lemma}
\begin{pf}
First, by (\ref{eq:k=1-L-const}) and the definitions of $\Delta_n$
and $\Delta_{n,\nu}$ in (\ref{eq:LL}) and (\ref{eq:LL-d}),
%
%
\begin{eqnarray} \label{eq:L-ratio-one}
\Delta_1(\rx)&\equiv&
\max_{a,b,c,d} \biggl|\frac{P_{01}(b,c)}{P_{01}(a,c)}- \frac
{P_{01}(b,d)}{P_{01}(a,d)}\biggr|\le\frac{\gamma(1-\phi_*)}{\phi_*},
\nonumber\\[-8pt]\\[-8pt]
\Delta_{1,\nu}(\rx)
&\equiv&0, \qquad\nu>0.\nonumber
\end{eqnarray}
By (\ref{eq:IL}), $I_{n,ec}^{(1)} = P_{n-1,n}(e,c)$, so
(\ref{eq:d-R}) gives $\Delta_n(\rx) \le\gamma\Delta_{n-1}(\rx)$.
Thus,
%
%
\begin{equation} \label{eq:Delta-gamma}
\Delta_n(\rx) \le\alpha\gamma^n\qquad \forall n\ge1,
\rx>0.
\end{equation}

Let $R_{n,\nu,c}(\rx)$ be as in (\ref{eq:d-lambda}) and
\[
\bar\Delta_{n,l}(\rx) = \max_{|\nu|=l} \Delta_{n,\nu}(\rx).
\]
Recall the definition of $V_n(\rx_0)$ in (\ref{eq:V-n-def}). For
brevity, write $v_n = V_n(\rx_0)$. By (\ref{eq:R-D-n}), there are
constants $c_l>1$, such that
%
%
\begin{equation}\label{eq:R-C-V}
{\max_{|\nu|=l, c}}
|R_{n,\nu,c}(\rx)| \le\half
c_l v_{n-1}^l \sum_{i=0}^{l-1}\bar\Delta_{n-1,i}(\rx),
\end{equation}
for $l=1,\ldots, q$, $n\ge1$, $\rx_0>0$ and $|\rx|\le\rx_0$. Then
by (\ref{eq:d-R-n}), for $n\ge0$,
%
%
\begin{equation} \label{eq:Delta-recur}
\bar\Delta_{n+1,l}(\rx) \le\gamma
\bar\Delta_{n,l}(\rx) + c_l v_n^l
\sum_{i=0}^{l-1} \bar\Delta_{n,i}(\rx).
\end{equation}

We show by induction that for $l\ge1$ and $n\ge0$,
%
%
\begin{equation} \label{eq:Delta-bound}
\bar\Delta_{n+1,l}(\rx) \le\alpha\gamma^{(n+1-l)\vee1}
n c_l v_n^l \prod_{i=1}^{l-1} (1+n c_i v_n^i),
\end{equation}
where $\bar\Delta_{n+1,0}(\rx) = \Delta_{n+1}(\rx)$.

By (\ref{eq:L-ratio-one}), (\ref{eq:Delta-bound}) holds for $n=0$
and $l\ge1$. Let $n\ge1$ next. If $l=1$, then by
(\ref{eq:Delta-gamma}) and (\ref{eq:Delta-recur}),
\[
\bar\Delta_{n+1,1}(\rx) \le\gamma\bar\Delta_{n,1}(\rx) +
c_1 v_n\Delta_n(\rx)
\le\gamma\bar\Delta_{n,1}(\rx) + \alpha\gamma^n c_1 v_n,
\]
and by induction on $n$ and (\ref{eq:L-ratio-one}),
\[
\bar\Delta_{n+1,1}(\rx)
\le
\gamma^n\bar\Delta_{1,1}(\rx) + \alpha\gamma^n c_1\sum_{s=1}^n v_s
=\alpha\gamma^n c_1\sum_{s=1}^n v_s
\le\alpha\gamma^n c_1 n v_n.
\]
So (\ref{eq:Delta-bound}) holds for $l=1$. Suppose
(\ref{eq:Delta-bound}) holds for $1\le l<k$. By $\gamma\in(0,1)$,
%
%
\begin{eqnarray}\label{eq:sum-D}
\sum_{i=0}^{k-1} \bar\Delta_{n,i}(\rx_0)
&
=&
\Delta_n(\rx_0)+ \sum_{i=1}^{k-1} \bar\Delta_{n,i}(\rx_0)
\nonumber
\\
&
\le&
\alpha\Biggl\{\gamma^n+ \sum_{i=1}^{k-1} \gamma^{(n-i)\vee1} c_i
(n-1)v_{n-1}^i \prod_{h=1}^{i-1} [1+c_h (n-1)v_{n-1}^h] \Biggr\}
\hspace*{-27pt}
\nonumber\\[-8pt]\\[-8pt]
&
\le&
\alpha\gamma^{(n+1-k)\vee1}
\Biggl\{1+ \sum_{i=1}^{k-1} c_i n v_n^i\prod_{h=1}^{i-1} (1+c_h n v_n^h) \Biggr\}
\nonumber
\\
&
=&
\alpha\gamma^{(n+1-k)\vee1} \prod_{i=1}^{k-1} (1+c_in v_n^i),\nonumber
\end{eqnarray}
so by (\ref{eq:Delta-recur}),
\[
\bar\Delta_{n+1,k}(\rx) \le\gamma\bar\Delta_{n,k}(\rx)
+ \alpha\gamma^{(n+1-k)\vee1} c_k v_n^k
\prod_{i=1}^{k-1} (1+ c_i n v_n^i).
\]
%
By induction on $n$, it is seen that $\bar{\Delta}_{n,k}(\rx)$ satisfies (\ref{eq:Delta-bound}).
As a by-product, by (\ref{eq:R-C-V}) and (\ref{eq:sum-D}),
%
%
\begin{equation} \label{eq:R-gamma}
{\max_{|\nu|=l, c}} |R_{n,\nu,c}(\rx)|\le
\half\alpha\gamma^{(n-l)\vee1} c_l v_{n-1}^l
\prod_{i=1}^{l-1} (1+c_i nv_{n-1}^i).
\end{equation}

Combining (\ref{eq:L-ratio-recurse}), (\ref{eq:Delta-bound}) and
(\ref{eq:R-gamma}), for any $|\nu|=l$,
\begin{eqnarray*}
&&\biggl|\biggl(\frac{L_{n,b c}}{L_{n,a c}}\biggr)^{(\nu)}(\rx) - \biggl(\frac{L_{n-1,b
c}}{L_{n-1,a c}}\biggr)^{(\nu)}(\rx)\biggr|\\
&&\qquad\le
\Delta_{n-1,l}(\rx) + 2|R_{n,\nu,c}(\rx)| \\
&&\qquad\le
\alpha\gamma^{(n-l-1)\vee1} n c_l v_{n-1}^l \prod_{i=1}^{l-1}
(1+c_i n v_{n-1}^i).
\end{eqnarray*}

Let $T_{n,\nu}$ be as in (\ref{eq:d-Lambda}). With
(\ref{eq:Delta-bound}) being established now, by (\ref{eq:T-D-n}),
we get the following bounds similar to (\ref{eq:R-gamma}):
%
%
\begin{equation} \label{eq:T-gamma}
{\max_{|\nu|=l}} |T_{n,\nu}(\rx)|\le\half
\alpha\gamma^{(n-l)\vee1} c_l v_n^l
\prod_{i=1}^{l-1} (1+n c_i v_{n-1}^i).
\end{equation}

Combine (\ref{eq:L-ratio-recurse}), (\ref{eq:d-Lambda-2}) and the
above inequalities. It is seen that for some constants $C_l>1$,
\[
\bigl|\Lambda_{n,a}^{(\nu)}- \Lambda_{n-1,a}^{(\nu)}\bigr|
\le C_l \alpha\gamma^{(n-l-1)\vee1} n^l v_n^{l(l+1)/2}.
\]
Then applying Lemma \ref{lemma:bound} to $v_n=V_n(\rx_0)$,
the lemma is proved.
\end{pf}

Now for $n\ge1$, $|\Lambda_{n,a}^{(\nu)}(\rx)| \le
|\Lambda_{1,a}^{(\nu)}(\rx)| + \sum_{k=2}^n
|\Lambda_{n,a}^{(\nu)}(\rx)-\Lambda_{n-1,a}^{(\nu)}(\rx)|$. Letting
$k=1$ in (\ref{eq:d-Lambda-2a}) and (\ref{eq:T-D-n}) and
combining them with (\ref{eq:k=1-L-const}) and
(\ref{eq:L-ratio-one}), it is seen that $|\Lambda_{1,a}^{(\nu)}(\rx
)|\le
C |V_1(\rx)|^{|\nu|}$, where $C$ is a constant. Together with
(\ref{eq:k=1-step}), this implies there is a constant $C_l =
C_l(\gamma, \phi_*)$, such that for $\nu$ with $1\le|\nu|=l\le q$,
%
%
\begin{equation} \label{eq:k=1-bound}
\bigl|\Lambda_{n,a}^{(\nu)}(\rx)\bigr|
\le C_l \sum_{t=1}^\infty\beta_{l,t}
[q+M_q(Z_t,\rx_0)]^{ql(l+1)/2}, \qquad|\rx|\le\rx_0,
\end{equation}
where $\beta_{l,t}= \sum_{k=t+1}^\infty\gamma^k k^{l(l+1)} =
o((c\gamma)^t)$ for any $0<c<1/\gamma$.

Part 2 of Theorem \ref{thm:k=1} is an immediate consequence of
the next result.
\begin{lemma} \label{lemma:Lambda-step2}
Let $\rx_0>0$. Almost surely, the following statements hold for all
$|\rx|\le\rx_0$, $n\ge1$ and $\nu$ with $1\le|\nu|\le q$.

1. $\mean[(\ln\Lambda_{n,a})^{(\nu)}(\rx)\mid\state]$ and $\mean
[(\ln
\Lambda_{n,a})(\rx)\mid\state]^{(\nu)}$ both exist and are equal.

2. As $n\toi$, $\mean[(\ln\Lambda_{n,a})^{(\nu)}(\rx)\mid\state]
\to
\mean[(\ln\sL_a)^{(\nu)}(\rx)\mid\state]$.

3. As $n\toi$, $(\mean[\ln\Lambda_{n,a}(\rx)\mid\state])^{(\nu)}
\to
(\mean[\ln\sL_a(\rx)\mid\state])^{(\nu)}$.
\end{lemma}
\begin{pf}
1. It is not hard to see that $(\ln\Lambda_{n,a})^{(\nu)}(\rx)$ is a
linear combination of products of the form
\[
h_{n,\nu_1, \ldots, \nu_s}(\rx):=\frac{\Lambda_{n,a}^{(\nu
_1)}(\rx)\cdots
\Lambda_{n,a}^{(\nu_s)}(\rx)}{\Lambda_{n,a}(\rx)^s},\qquad
\nu_k> 0,\qquad \nu_1+\cdots+\nu_s=\nu.
\]
By (\ref{eq:Lambda-ul}) and (\ref{eq:k=1-bound}),
\[
|h_{n,\nu_1, \ldots, \nu_s}(\rx)|\le
\zeta:
=C\prod_{k=1}^s \sum_{t=1}^\infty\beta_{l,t}
[q+M_q(Z_t,\rx_0)]^{q|\nu_k|(|\nu_k|+1)/2},\qquad |\rx|\le\rx_0,
\]
with $C = C(\gamma,\phi_*)$ a constant. As $\sum_k
|\nu_k|(|\nu_k|+1) \le|\nu|(|\nu|+1)$, by
Assumption~\hyperlink{A:density2}{5} and the independence of $Z_t$,
$\mean\zeta<\infty$. Note that $\zeta$ is independent of
$\state$. It follows that $(\ln\Lambda_{n,a})^{(\nu)}(\rx)$ for all
$n$ and $|\rx|\le\rx_0$ are bounded by a single random variable
that has a finite expectation and is independent of $\state$. This
implies $\mean[(\ln\Lambda_{n,a})^{(\nu)}(\rx)\mid\state]$ exists, and
together with $\ln\Lambda_{n,a}\in C^{(q)}$, implies the rest of part
1 through dominated convergence.

2. By Theorems \ref{thm:LR-limit} and \ref{thm:L-limit},
$\Lambda_{n,a}^{(\nu)}(\rx)$ converges as $n\toi$ for all $\rx$. By
Lem\-ma~\ref{lemma:Lambda-step} and (\ref{eq:Lambda-ul}),
$(\ln\Lambda_{n,a})^{(\nu)}(\rx)$ converges. Then the claim follows
from dominated convergence.

3. Consider $h_{n,\nu_1, \ldots, \nu_s}(\rx)$ again. By Lemma
\ref{lemma:Lambda-step} and (\ref{eq:Lambda-ul}), it can be seen
that for $\nu_1, \ldots, \nu_s>0$ with $\nu_1+\cdots+\nu_s=\nu$,
$|h_{n+1,\nu_1, \ldots, \nu_s}(\rx)-h_{n,\nu_1, \ldots, \nu_s}(\rx
)|\le C \gamma_1^n \zeta$
holds for $|\rx|\le\rx_0$, where $C>0$ $\gamma_1\in(\gamma,1)$ are constants and $\zeta>0$
is a random variable independent of $\state$ with
$\mean\zeta<\infty$. As a result,
$\mean[(\ln\Lambda_{n,a})^{(\nu)}(\rx)\mid\state]$ converges uniformly
on each compact set of $\rx$. Together with part 1, this implies
part~3.
\end{pf}

\subsection{Proof for the binary case}
The following simple identity will be repeatedly used. For any
function $F$ on $\{0,1\}$, denote $dF=F(1)-F(0)$. Then for $s$,
$t\in\Ints$,
%
%
\begin{eqnarray}
\label{eq:contrast}
\mean_\stx[F(\stx_t)\mid\stx_s=1] -
\mean_\stx[F(\stx_t)\mid\stx_s=0]
&=&D_{st} \,d F, \\
\label{eq:contrast2}
F(0) - \mean_\stx[F(\stx_t)\mid\stx_s=0] &=& -P_{st}(0,1) \,d F.
\end{eqnarray}

Define for $t\in\Ints$ and $n\ge1$,
\[
\ell_t(\rx, a) = \ln\psi_t(\rx,a),\qquad
S_n(\rx)=\sum_{t=1}^n \ell_t(\rx,\stx_t).
\]
Then $\lambda_n(\rx)=\ln\mean_\stx[e^{S_n(\rx)}\mid\stx_0=1] -\ln
\mean_\stx[e^{S_n(\rx)}\mid\stx_0=0]$ in (\ref{eq:lambda-n}).

\begin{pf*}{Proof of Theorem \ref{thm:d-L-1-2}}
For $n\ge1$, by (\ref{eq:contrast}),
\begin{eqnarray*}
\lambda_n'(0)
&=&
\mean_\stx[S_n'(0)\mid\stx_0=1] -
\mean_\stx[S_n'(0)\mid\stx_0=0]
\\
&=&
\sum_{t=1}^n \{\mean_\stx[\ell_t'(0,\stx_t)\mid\stx_0=1] - \mean
_\stx[\ell_t'(0,\stx_t)\mid\stx_0=0] \} \\
&=& \sum_{t=1}^n D_{0t} d_t'(0).
\end{eqnarray*}
By Theorems \ref{thm:L-limit} and \ref{thm:d-L-limit}, letting
$n\toi$, (\ref{eq:d-L-1}) follows. To get $\sr''(0)$, for $n\ge1$,
\begin{eqnarray*}
\lambda_n''(0)
&=&\mean_\stx[S_n''(0)\mid\stx_0=1] -
\mean_\stx[S_n''(0)\mid\stx_0=0]
+ \var_\stx[S_n'(0)\mid\stx_0=1] \\
&&{} - \var_\stx[S_n'(0)\mid\stx_0=0].
\end{eqnarray*}
Similar to the calculation of $\sr'(0)$,
\[
\lim_{n\toi}
\{\mean_\stx[S_n''(0)\mid\stx_0=1] - \mean_\stx[S_n''(0)\mid\stx
_0=0] \}=
\sum_{t=1}^\infty D_{0t} d_t''(0).
\]

On the other hand, denoting by $\delta_t$ the random variable
$\ell_t'(0,\stx_t)$,
\[
\var_\stx[S_n'(0)\mid\stx_0]
= \sum_{t=1}^n \var_\stx(\delta_t\mid\stx_0)
+ 2\sum_{1\le s<t\le n}
\cov_\stx(\delta_s, \delta_t\mid\stx_0).
\]
Given $1\le s\le t\le n$, let $F(\stx_s) = \delta_s
\mean_\stx[\delta_t\mid\stx_s]$. By $\mean_\stx[\delta_s
\delta_t\mid\stx_0] =\mean_\stx[F(\stx_s)\mid\stx_0]$ and
(\ref{eq:contrast}),
\[
\mean_\stx[\delta_s\delta_t\mid\stx_0=1]-
\mean_\stx[\delta_s\delta_t\mid\stx_0=0]= D_{0s} \,d F.
\]
Similarly, by (\ref{eq:contrast}), $\mean_\stx[\delta_t\mid\stx_s=1] =
\mean_\stx[\delta_t\mid\stx_s=0] + D_{st} d_t'(0)$.
Then, as $\ell_s'(0,1) = \ell_s'(0,0)+d_s'(0)$,
\begin{eqnarray*}
d F
&=& F(1) - F(0)
= \ell_s'(0,1) \mean_\stx(\delta_t\mid\stx_s=1) -
\ell_s'(0,0) \mean_\stx(\delta_t\mid\stx_s=0) \\
&=& \mean_\stx(\delta_t\mid\stx_s=0)d_s'(0)
+ D_{st} \ell_s'(0,0)d_t'(0)
+ D_{st} d_s'(0) d_t'(0)
\end{eqnarray*}
and likewise,
\begin{eqnarray*}
&&\mean_\stx(\delta_s\mid\stx_0=1)
\mean_\stx(\delta_t\mid\stx_0=1)
-
\mean_\stx(\delta_s\mid\stx_0=0)
\mean_\stx(\delta_t\mid\stx_0=0) \\
&&\qquad=
D_{0s} \mean_\stx(\delta_t\mid\stx_0=0) d_s'(0)
+ D_{0t} \mean_\stx(\delta_s\mid\stx_0=0) d_t'(0)
+ D_{0s} D_{0t} d_s'(0) d_t'(0).
\end{eqnarray*}

Combining the above identities,
\[
\cov_\stx(\delta_s, \delta_t\mid\stx_0=1)
-\cov_\stx(\delta_s, \delta_t\mid\stx_0=0) = I_1 + I_2 + I_3, \\
\]
with
\[
\cases{
I_1 =
D_{0s}[\mean_\stx(\delta_t\mid\stx_s=0) - \mean_\stx(\delta_t\mid\stx
_0=0) ] d_s'(0), \cr
I_2 =
[D_{0s}D_{st} \ell_s'(0,0) - D_{0t}\mean_\stx(\delta_s\mid\stx_0=0)
] d_t'(0), \cr
I_3 =
D_{0s}(D_{st} - D_{0t}) d_s'(0) d_t'(0).}
\]
By conditioning on $\stx_s$,
\begin{eqnarray*}
&&\mean_\stx(\delta_t\mid\stx_s=0) -
\mean_\stx(\delta_t\mid\stx_0=0)\\
&&\qquad=
\mean_\stx(\delta_t\mid\stx_s=0) -
\mean_\stx[\mean_\stx(\delta_t\mid\stx_s)\mid\stx_0=0] \\
&&\qquad\stackrel{(\mathrm{a})}{=}
-P_{0s}(0,1) [\mean_\stx(\delta_t\mid\stx_s=1)
-\mean_\stx(\delta_t\mid\stx_s=0)]\\
&&\qquad\stackrel{(\mathrm{b})}{=}
- D_{st} P_{0s}(0,1) d_t'(0),
\end{eqnarray*}
where (a) is due to (\ref{eq:contrast2}), and (b) due to
(\ref{eq:contrast}). By (\ref{eq:Dproduct}), $D_{0s} D_{st} =
D_{0t}$. Therefore, $I_1 = - D_{0t} P_{0s}(0,1) d_s'(0) d_t'(0)$.
Likewise,
\[
I_2 = D_{0t}[\ell_s'(0,0) - \mean_\stx(\delta_s\mid\stx_0=0)]
d_t'(0) = -D_{0t} P_{0s}(0,1) d_s'(0) d_t'(0)
\]
and $I_3 = D_{0t}(1-D_{0s})d_s'(0)d_t'(0)$. Then
(\ref{eq:d-L-2}) follows from
\begin{eqnarray*}
&&\cov_\stx(\delta_s, \delta_t\mid\stx_0=1)
-\cov_\stx(\delta_s, \delta_t\mid\stx_0=0)\\
&&\qquad= D_{0t}[P_{0s}(1,0) - P_{0s}(0,1)] d_s'(0) d_t'(0)
\end{eqnarray*}
and Theorems \ref{thm:L-limit} and \ref{thm:d-L-limit}.
\end{pf*}

To prove the rest of the results, recall $\lambda(x,\px)=\ln
f(x,\px)$.
\begin{pf*}{Proof of Proposition \ref{prop:d-uv}}
Given $t$, $Z$ and $\state$, $\ell_t(\rx, a)$ is a composite of
functions $\lambda(x,\px)$, $\tf(Z_t, v)$, $\prx_a(\rx)$ and
$\prx_{\state_t}(\rx)$, such that
\[
\ell_t(\rx,a) = \lambda(\tf(Z_t, \prx_{\state_t}(\rx)),
\prx_a(\rx)),
\]
so by the chain rule for differentiation,
\[
\ell_t'(\rx,a)
=
\frac{\partial\lambda(x,\px)}{\partial x}\,
\frac{\partial\tf(Z_t,v)}{\partial v}\prx_{\state_t}'(\rx)
+
\frac{\partial\lambda(x,\px)}{\partial\px}\prx_a'(\rx),
\]
where the right-hand side is evaluated at $x = \tf(Z_t, v)$,
$v=\prx_{\state_t}(\rx)$, and $\px=\prx_a(\rx)$. Since $\prx
_1(0) =
\prx_0(0)=0$, the first summand on the right-hand side takes the
same value for $a=0$, 1. Therefore, (\ref{eq:uv-d-diff-1}) holds.

Likewise,
\begin{eqnarray*}
\ell_t''(\rx,a)
&=&
\frac{\partial^2 \lambda}{\partial x^2}
\biggl[\frac{\partial\tf}{\partial v}\biggr]^2
\prx_{\state_t}'(\rx)^2
+ 2\,\frac{\partial^2 \lambda}{\partial x\,\partial\px}\,
\frac{\partial\tf}{\partial v}
\prx_{\state_t}'(\rx)\prx_a'(\rx)
+\frac{\partial^2 \lambda}{\partial\px^2}\prx_a'(\rx)^2 \\
&&{}+\frac{\partial\lambda}{\partial x}\,
\frac{\partial^2\tf}{\partial v^2}\prx_{\state_t}'(\rx)^2
+\frac{\partial\lambda}{\partial x}\,
\frac{\partial\tf}{\partial v}\prx_{\state_t}''(\rx)
+ \frac{\partial\lambda}{\partial\px}\prx_a''(\rx),
\end{eqnarray*}
where again the right-hand side is evaluated at $x = \tf(Z_t, v)$,
$v=\prx_{\state_t}(\rx)$, and $\px=\prx_a(\rx)$. Then
(\ref{eq:uv-d-diff-2}) follows.
\end{pf*}

\begin{pf*}{Proof of Proposition \ref{prop:E-d-uv}}
We shall first show for any $t$,
%
%
\begin{eqnarray}
\label{eq:E-diff-1}
\mean[d_t'(0)\mid\state] &=& 0, \\
\label{eq:Var-diff-1}
\var[d_t'(0)\mid\state] &=& [\prx_1'(0)-\prx_0'(0)]^2 J(0),
\\
\label{eq:E-diff-2}
\mean[d_t''(0)\mid\state]
&=& [\prx_1'(0)-\prx_0'(0)] [2\prx_{\state_t}'(0)-\prx_0'(0) -
\prx_1'(0)] J(0).
\end{eqnarray}

Denote $\xi_t = \tf(Z_t,0)$. Then $\xi_t$ has density $f(x,0)$ and
log-density $\lambda(x,0)$. Take expectation conditional on
$\state$ on both sides of (\ref{eq:uv-d-diff-1}) to get
\[
\mean[d_t'(0)\mid\state] = [\theta_1'(0) - \theta_0'(0)]
\mean\biggl[\frac{\partial\lambda(\xi_t,0)}{\partial\px}\biggr].
\]
Then (\ref{eq:E-diff-1}) follows from the property of score
function.

For the same reason, (\ref{eq:Var-diff-1}) follows as well and,
taking expectation conditional on $\state$ on both sides of
(\ref{eq:uv-d-diff-2}),
\begin{eqnarray*}
\mean[d_t''(0)\mid\state]
&=& 2[\prx_1'(0)-\prx_0'(0)] \prx_{\state_t}'(0)
\mean\biggl[\frac{\partial^2 \lambda(\xi_t,0)}{\partial
x\,\partial\px}\,
\frac{\partial\tf(Z_t,0)}{\partial v} \biggr]\\
&&{}-[\prx_1'(0)^2-\prx_0'(0)^2] J(0).
\end{eqnarray*}

Therefore, to prove (\ref{eq:E-diff-2}), it suffices to show
%
%
\begin{equation} \label{eq:FI}
\mean\biggl[\frac{\partial^2 \lambda(\xi_t,0)}{\partial x\,\partial\px}\,
\frac{\partial\tf(Z_t,0)}{\partial v} \biggr] = J(0).
\end{equation}
Define
\[
g(v,Z_t)
=
\frac{\partial\lambda(\tf(Z_t, v), \px)}{\partial\px}
\bigg|_{\px=0}
= \frac{1}{f(\tf(Z_t, v), 0)}\,
\frac{\partial f(\tf(Z_t,v),0)}{\partial\px}.
\]
Observe that
\[
\frac{\partial g(v, Z_t)}{\partial v}\bigg|_{v=0}=
\frac{\partial^2 \lambda(\xi_t,0)}{\partial x\,\partial\px}\,
\frac{\partial\tf(Z_t,0)}{\partial v}.
\]
Therefore, the left-hand side of (\ref{eq:FI}) is equal to
\[
\mean\biggl[\frac{\partial g(v,Z_t)}{\partial v} \bigg|_{v=0} \biggr]
=
\frac{\partial\mean[g(v,Z_t)]}{\partial v}\bigg|_{v=0}.
\]

By assumption \hyperlink{A:noise}{1}, $\tf(Z_t,v)$ has density $f(x,v)$.
Therefore, the right-hand side of the above identity is equal to
\[
\frac{\partial}{\partial v}
\biggl[\int\frac{1}{f(x, 0)}\,\frac{\partial f(x,0)}{\partial\px} f(x,v) \,d
x \biggr]_{v=0}
=
\int\frac{1}{f(x, 0)}\biggl[\frac{\partial
f(x,0)}{\partial\px}\biggr]^2\,
d x = J(0),
\]
which gives (\ref{eq:FI}).

From Theorem \ref{thm:k=1}, (\ref{eq:d-L-1}) and
(\ref{eq:E-diff-1}),
\[
\mean[\sr'(0)\mid\state] = \sum_{t=1}^\infty D_{0t}\mean[d_t'(0)\mid
\state] = 0
\]
showing (\ref{eq:E-d-L-1}). On the other hand, given $\state$,
since $Z_t$ are independent, $d_s'(0)$ are independent of $d_t'(0)$
for $s<t$. Then by $\mean[d_t'(0)\mid\state]=0$ and
(\ref{eq:d-L-2}),
\begin{eqnarray*}
\mean[\sr''(0)\mid\state]
&=&\sum_{t=1}^\infty D_{0t}\{\mean[d_t''(0)\mid\state] +
[P_{0t}(1,0)-P_{0t}(0,1)] \var[d_t'(0)\mid\state] \} \\
&=&[\prx_1'(0)-\prx_0'(0)] J(0)
\sum_{t=1}^\infty D_{0t} f_t,
\end{eqnarray*}
where
\begin{eqnarray*}
f_t
&=&2\prx_{\state_t}'(0) - \prx_0'(0) - \prx_1'(0)
+
[P_{0t}(1,0)-P_{0t}(0,1)] [\prx_1'(0) - \prx_0'(0)] \\
&=&
[\prx_1'(0)-\prx_0'(0)][2\state_t-P_{0t}(1,1)-P_{0t}(0,1)].
\end{eqnarray*}
Therefore, (\ref{eq:E-d-L-2}) holds.
\end{pf*}

\vspace*{-15pt}

\begin{appendix}\label{app}
\section*{Appendix}
In this Appendix, we make a general statement on the conditional
likelihood under the FDR criterion. Let $H_1, \ldots, H_{n}$ be a set of
hypotheses being tested and let $X$ be the available data. Let
$p_k = \mathsf{Pr}\{H_k \mbox{ is false}\mid X\}$. For any testing procedure
based on $X$, let $R$ be the total number of rejected $H_k$ and $V$
that of rejected true~$H_k$. Then the number of rejected false nulls
is $R-V$.
\setcounter{theorem}{0}
\begin{prop}\label{prop:FDR}
Given $\alpha\in(0,1)$, among all testing procedures whose
rejection decisions are uniquely determined by $X$ and
which satisfy the FDR control criterion
\[
\mathrm{FDR}=\mean\biggl[\frac{V}{R\vee1}\Big| X\biggr] \le\alpha,
\]
the following Benjamini--Hochberg procedure
\cite{benjaminihoc95} has the largest $\mean[R-V\mid X]$:
\begin{enumerate}
\item sort $q_i=1-p_i$ into $q_{(1)}\le q_{(2)}
\le\cdots\le q_{(n)}$;
\item let $r=\max\{j\dvtx q_{(1)}+\cdots+q_{(j)} \le\alpha j\}$ and
reject $H_k$ if $q_k\le q_{(r)}$.
\end{enumerate}
\end{prop}
\begin{pf}
Given a procedure with $R>0$, let $H_{i_k}$, $k=1,\ldots, R$ be
the rejected nulls. Then, as in \cite{chi08mvp},
$\mathrm{FDR} = \sum_{j=1}^R q_{i_j}/R\ge\sum_{j=1}^R q_{(j)}/R$,
while $\mean[R-V\mid X] = R - \sum_{j=1}^R q_{i_j}\le R -
\sum_{j=1}^R q_{(j)}$. It is then not hard to see that under the FDR
control criterion, the procedure in the proposition attains the
largest value of $E[R-V\mid X]$.
\end{pf}
\end{appendix}


%
\printaddresses

\end{document}